\documentclass[10pt]{amsart}

\usepackage{amsfonts, amssymb, amsxtra, amstext, amsbsy, amscd, amsmath, commath, eucal, graphicx, hyperref, latexsym, mathrsfs, mathtools, tikz-cd, xcolor}
\usepackage[utf8]{inputenc}
\usepackage[all]{xy}

\newtheorem{theorem}{Theorem}[section]
\newtheorem{lemma}[theorem]{Lemma}
\newtheorem{proposition}[theorem]{Proposition}
\newtheorem{corollary}[theorem]{Corollary}
\theoremstyle{definition}

\newtheorem{example}[theorem]{Example}
\newtheorem{remark}[theorem]{Remark}

\numberwithin{equation}{section}
\newtheorem{question}[theorem]{Question}

\def\C{\mathop{\mathsf{C}\hspace{0mm}}\nolimits}
\def\cf{\mathop{\rm cf}\nolimits}

\def\CN{\mathop{\mathsf{CN}}\nolimits}

\def\img{\mathop{\rm img}\nolimits}
\def\int{\mathop{\rm int}\nolimits}

\def\P{\mathop{\mathsf{P}\hspace{0mm}}\nolimits}
\def\Pot{\mathop{\mathscr{P}\hspace{0mm}}\nolimits}
\def\R{\mathop{\mathsf{R}\hspace{0mm}}\nolimits}

\def\WP{\mathop{\mathsf{WP}}\nolimits}

\newcommand\restr[2]{\ensuremath{#1\!\!\upharpoonright_{#2}}}

\usepackage{lineno}

\let\mathcal\mathscr

\makeatletter
\@namedef{subjclassname@2020}{\textup{2020} Mathematics Subject Classification}
\makeatother

\begin{document}


\title{Some notes on topological calibers}

\author{Alejandro R\'ios-Herrej\'on}
\address {A. Ríos Herrej\'{o}n\\
Departamento de Matem\'aticas, Facultad de Ciencias, Universidad Nacional Aut\'onoma de M\'exico, Circuito ext. s/n, Ciudad Universitaria, C.P. 04510,  M\'exico, CDMX}
\email{chanchito@ciencias.unam.mx}

\thanks{The first author received grant 814282 from CONACYT}

\author[\'A. Tamariz-Mascar\'ua]{\'Angel Tamariz-Mascar\'ua}
\address{Á. Tamariz Mascar\'ua, Departamento de Matem\'aticas, Facultad de Ciencias, Universidad Nacional Aut\'onoma de M\'exico, Circuito ext. s/n, Ciudad Universitaria, C.P. 04510,  M\'exico, CDMX}
\email{atamariz@unam.mx}
\urladdr{}



\begin{abstract}
We show that the definition of caliber given by Engelking in  \cite{engelking1989}, which we will call {\it caliber*}, differs from the traditional notion of this concept in some cases and agrees in others.

For instance, we show that if $\kappa$ is an infinite cardinal with $2^{\kappa}<\aleph_\kappa$ and $\cf(\kappa)>\omega$, then there exists a compact Hausdorff space $X$ such that $o(X)=2^{\aleph_\kappa}=|X|$, $\aleph_\kappa$ is a caliber* for $X$ and $\aleph_\kappa$ is not a caliber for $X$.

On the other hand, we obtain that if $\lambda$ is an infinite cardinal number, $X$ is a Hausdorff space with $|X|>1$, $\phi\in \{w ,nw\}$, $o(X) = 2^{\phi(X)}$ and $\mu := o\left(X^\lambda\right)$, then the calibers of $X^\lambda$ and the true calibers* (that is, those which are less than or equal to $\mu$) coincide, and are precisely those that have uncountable cofinality.
\end{abstract}

\maketitle

\section{Introduction}

The classical definition of caliber for a topological space states that if $\kappa$ is an infinite cardinal, then $\kappa$ is a caliber for $X$ if for every collection $\{U_\alpha : \alpha<\kappa\}$ of non-empty open subsets of $X$, there exists $J\subseteq \kappa$ such that $\abs{J}=\kappa$ and $\{U_\alpha : \alpha\in J\}$ satisfies $\bigcap_{\alpha\in J} U_\alpha \neq \emptyset$ (see \cite{sanin1948}).

However, despite the fact that this convention seems to be the most used (see, for example, \cite{comneg1982}, \cite{juhasz1980} and \cite{shelah1977}), in \cite[2.7. 11, p.~116]{engelking1989} Engelking defines: an infinite cardinal $\kappa$ is a caliber for a topological space $X$ if for every collection $\mathcal{U}$ of cardinality $\kappa$ of non-empty open subsets of $X$, there exists a subcollection $\mathcal{V}\subseteq \mathcal{U}$, such that $\abs{\mathcal{V}}=\kappa$ and $\bigcap\mathcal{V}\neq\emptyset$.

Naturally, an immediate observation is that in \v{S}anin's definition repetitions in the enumeration are allowed; therefore, this is an indicator that these two notions of caliber may not coincide. The main goal of this paper is to delve into the differences that lie between these concepts and to expose that they indeed differ in some aspects.

We will call an infinite cardinal number $\kappa$ {\it caliber*} if it satisfies the definition given by Engelking, and we will call it simply {\it caliber} if it satisfies \v{S}anin's definition.

In Section 2 we list, without proof, some basic known results on calibers. Section 3 is devoted to obtaining basic results on calibers* and some of their variations such as precalibers* and weak  precalibers*. Furthermore, in this section, we present examples of topological spaces for which their non trivial calibers* (those that do not exceed the number of open subsets) do not coincide with their calibers. Finally, we calculate calibers* in concrete classes of topological spaces like those  hyperconnected and in $\beta \omega$. In Section 4, we study the calibers* for topological sums, and in Section 5 we analyze calibers* for topological products.

\section{Preliminaries}

Any topological or set-theoretic concept that is not explicitly defined in this text should be understood as in \cite{engelking1989} and \cite{kunen1980}, respectively.

For a topological space $X$, the symbol $\tau_{X}$ is the topology of $X$, and we will denote by $\tau_{X}^{+}$ the set $\tau_{X}\setminus \{\emptyset\}$.
For a cardinal number $\kappa$, we will use the symbol $D(\kappa)$ to denote the discrete space of cardinality $\kappa$.

If $\alpha$ is an ordinal number, $\aleph_\alpha$ will represent the cardinal number of the set $\omega_\alpha$ (see \cite{kunen1980}). Furthermore, we will use the symbol 
$\omega$ to denote the first infinite ordinal and, therefore, the first infinite cardinal. The cardinal number of $\mathbb{R}$, the set of real numbers, will be denoted by $\mathfrak{c}$.

If $X$ is a set, the expression \lq\lq $X$ is a countable set{\rq\rq} will mean that there exists an injective function from $X$ into $\omega$. If $\kappa$ is a cardinal number, we will use the symbol $[X]^{<\kappa}$ to refer to the collection of subsets of $X$ of size (cardinality) less than $\kappa$. In the same way, $[X]^{\leq\kappa}$ will represent the family of all subsets of size at most $\kappa$, and we will denote by $[X]^{\kappa}$ the set $[X]^{\leq\kappa} \setminus [X]^{<\kappa}$. Finally, $\Pot(X)$ will be the power set of $X$.

The \textit{cofinality} of a cardinal $\kappa$, $\cf(\kappa)$, is the minimum ordinal number for which there exists a function $f:\cf(\kappa) \to \kappa$ so that for every $\alpha<\kappa$, there exists $\beta<\cf(\kappa)$ with $\alpha<f(\beta)$. Naturally, we always
have $\cf(\kappa)\leq \kappa$. We will say that an infinite cardinal $\kappa$ is \textit{regular} if $\kappa=\cf(\kappa)$; otherwise, we say that $\kappa$ is \textit{singular}. Lastly, we will denote by $\kappa^{+}$ the successor cardinal of $\kappa$; that is, the first cardinal $\lambda$ that satisfies the relation $\kappa<\lambda$.

Throughout this work we will denote by $\CN$ the proper class formed by the infinite cardinal numbers. For a topological space $X$, if $\mathcal{U}$ is a pairwise disjoint collection of $\tau_{X}^{+}$, we say that $\mathcal{U}$ is a \textit{cellular family} in $X$.
The terms {\it centered family} and {\it family with the finite intersection property} in $X$ mean the same; that is, they are used to designate a collection of subsets such that the intersection of every non-empty finite subcollection is non-empty. We say that a subset $\mathcal{U}$ of $\tau_X^+$ is {\it linked} if we have $ U\cap V\neq \emptyset$ for any $U,V\in \mathcal{U}$.

A cardinal number $\kappa$ is a {\it caliber} (resp., {\it precaliber}; {\it weak precaliber}) of a topological space $X$ if for every family $\{U_\alpha : \alpha < \kappa\}$ of non-empty open subsets of $X$, there exists $J \subseteq \kappa$ such that $|J| = \kappa$ and $\{U_\alpha : \alpha \in J \}$ satisfies 
$\bigcap_{\alpha \in J}U_\alpha \not= \emptyset$ (resp., it is centered; it is linked).

We will be working throughout this paper with the following collections of cardinal numbers: \begin{align*} \WP(X)&:=\{\kappa\in\CN : \kappa \ \text{is a weak precaliber for}\ X\}; \\
\P(X)&:=\{\kappa\in\CN : \kappa \ \text{is a precaliber for} \ X\}; \ \text{and} \\
\C(X)&:=\{\kappa\in\CN : \kappa \ \text{is a caliber for} \ X\}.
\end{align*}

Our basic reference texts for material related to topological cardinal functions will be \cite{hodel1984} and \cite{juhasz1980}. However, for the purposes of this text, the symbols $o(X)$ and $|X|$ will represent, correspondingly, the cardinality of $\tau_X$ and the cardinality of $X$; that is, we will dispense adding $\omega$ to these two cardinal numbers.

We compile below some known results that will be useful in the development of this text.

\begin{lemma}\label{lema_cofinalidad} If $X$ is a topological space and $\kappa$ is a caliber (resp., precaliber; weak precaliber) for $X$, then so is $\cf(\kappa)$.
\end{lemma}

\begin{lemma}\label{lema_densidad} 
If $X$ is a topological space, $\kappa$ is a cardinal number, and $d(X)<\cf(\kappa)$, then $\kappa$ is a caliber for $X$. Consequently, $\{\kappa\in\CN : \cf(\kappa)>d(X)\}\subseteq \C(X)$.
\end{lemma}

\begin{lemma}\label{precalibre_debil_celularidad} 
If $X$ is a topological space, $\kappa$ is a weak precaliber for $X$, and $\mathcal{U}$ is a cellular family in $X$, then $|\mathcal{U}|<\kappa $.
\end{lemma}

\begin{proposition}\label{calibre_subconjunto_denso} 
If $X$ is a topological space, $D$ is a dense subset of $X$, and $\kappa$ is an infinite cardinal, then the following statements are true.
\begin{enumerate}
\item $\kappa$ is a precaliber (resp., weak precaliber) for $D$ if and only if $\kappa$ is a precaliber (resp., weak precaliber) for $X$.
\item If $\kappa$ is a caliber for $D$, then $\kappa$ is a caliber for $X$.
\end{enumerate}
\end{proposition}

\begin{proposition}\label{calibre_subconjunto_abierto} 
If $X$ is a topological space, $U$ is an open subset of $X$, and $\kappa$ is a caliber (resp., precaliber; weak precaliber) for $X$, then $\kappa$ is a caliber (resp., precaliber; weak precaliber) for $U$.
\end{proposition}

\begin{lemma}\label{lema_funciones} 
Let $f:X\to Y$ be a continuous and surjective function. If $\kappa$ is a caliber (resp., precaliber, weak precaliber) for $X$, then $\kappa$ is a caliber (resp., precaliber; weak precaliber) for $Y$.
\end{lemma}

\begin{lemma}\label{lema_cubo_calibre} 
If $\kappa$ and $\lambda$ are cardinal numbers with $\cf(\kappa)>\omega$, then $D(2)^{\lambda}$ has caliber $\kappa$ (see \cite{shelah1977}).
\end{lemma}

\begin{lemma}\label{lema_cubo_topologia} 
For any cardinal number $\lambda$, $o\left(D(2)^{\lambda}\right)=2^{\lambda}$.
\end{lemma}

\section{Chain* conditions}

A cardinal number $\kappa$ is a {\it caliber*} (resp., {\it precaliber*}; {\it weak precaliber*}) of the topological space $X$ if for every family $\mathcal{U} $ of cardinality $\kappa$ of non-empty open subsets of $X$, there exists a subcollection $\mathcal{V} \subseteq \mathcal{U}$, such that $|\mathcal{V}| = \kappa$ and $\mathcal{V}$ has non-empty intersection (resp., is centered; is linked). In the remainder of this text we will be working with the following collections:

\begin{align*} \WP^*(X)&:=\{\kappa\in\CN : \kappa \ \text{is a weak precaliber* for} \ X\}; \\
\P^*(X)&:=\{\kappa\in\CN : \kappa \ \text{is a precaliber* for} \ X\}; \ \text{and} \\ 
\C^*(X)&:=\{\kappa\in\CN : \kappa \ \text{is a caliber* for} \ X\}.
\end{align*}

In this section we will present several basic results regarding the relations between the chain conditions that we have presented. Let us start with the following result.

\begin{lemma}\label{equiv_calibre_indices} 
Let $X$ be a topological space and $\kappa$ a regular cardinal. Then $\kappa$ is a caliber (resp., precaliber; weak precaliber) for $X$ if and only if $\kappa$ is a caliber* (resp., precaliber*; weak precaliber*).
\end{lemma}

\begin{proof} We give a proof for the case of caliber.
The other two cases can be similarly proved. Note that the direct implication is immediate and does not need the regularity of $\kappa$; it is sufficient to enumerate a $\mathcal{U}\in[\tau_{X}^{+}]^{\kappa}$ without repetitions. To verify the converse implication, suppose that $\{U_\alpha : \alpha<\kappa\}$ is a subset of $\tau_X^{+}$ and define $f:\kappa\to\{U_\alpha : \alpha<\kappa\}$ by $f(\alpha):=U_\alpha$. If $\left| \{U_\alpha : \alpha<\kappa\} \right| <\kappa$, then the regularity of $\kappa$ guarantees that there exists $\beta <\kappa$ such that $\left|f^{-1}\{U_\beta\}\right|=\kappa$. The collection $\{U_\alpha : \alpha\in f^{-1}\{U_\beta\}\}=\{U_\beta\}$ satisfies that $\bigcap\{U_\alpha : \alpha \in f^{-1}\{U_\beta\}\}\neq\emptyset$. Now if $\left| \{U_\alpha : \alpha<\kappa\} \right| =\kappa$, then the hypothesis gives us $\mathcal{V}\in [\{U_\alpha : \alpha<\kappa\}]^{\kappa}$ such that $\bigcap\mathcal{V}
\neq\emptyset$. Finally, if $J:=\{\alpha< \kappa : U_\alpha \in \mathcal{V}\}$, then $J\in [\kappa]^{\kappa}$ because the function $g :J\to\mathcal{V}$ given by $g(\alpha):=U_\alpha$ is surjective, and $\bigcap\{U_\alpha : \alpha\in J\}\neq\emptyset$.
\end{proof}

An equivalence of the chain conditions in terms of indexed families will help.

\begin{lemma}\label{lema_equivalencia*}  If $X$ is a topological space and $\kappa$ is a cardinal number, then the following statements are equivalent.

\begin{enumerate}
\item $\kappa$ is a caliber* (resp., precaliber*; weak precaliber*) for $X$.

\item For any family $\{U_\alpha : \alpha<\kappa\}\subseteq \tau_X^+$ enumerated without repetitions there exists $J\in [\kappa]^{\kappa}$ such that $\{U_\alpha : \alpha\in J\}$ has non-empty intersection (resp., is centered; is linked).

\item If $\{U_\alpha : \alpha<\kappa\}$ is a subset of $\tau_X$ with $\left|\{U_\alpha : \alpha<\kappa\}\right|=\kappa $, then there exists $J\in [\kappa]^{\kappa}$ with the following characteristics:
\begin{enumerate}
\item $\{U_\alpha : \alpha\in J\}=\{U_\alpha : \alpha<\kappa\}\setminus\{\emptyset\}$;
\item if $\alpha,\beta\in J$ are different, then $U_\alpha \neq U_\beta$; and
\item $\{U_\alpha : \alpha\in J\}$ has a non-empty intersection (resp., is centered; is linked).
\end{enumerate}
\end{enumerate}
\end{lemma}

\begin{proof} 
Since the implications $(1)\to (2)$ and $(3)\to (1)$ are evident, we only need to argue $(2) \to (3)$. Let us set $\{U_\alpha : \alpha<\kappa\}\subseteq\tau_X$ with $\left|\{U_\alpha : \alpha<\kappa\}\right|=\kappa$ and define an equivalence relation $\sim$ on $\kappa$ as follows: $\alpha\sim\beta$ if and only if $U_\alpha = U_\beta$. Then, if $e: \kappa/\!\!\sim \,\, \to \bigcup \kappa$ is a choice function and $I:= \img(e)$, then it is easy to verify that $I\in [\kappa]^{ \kappa}$, $\{U_\alpha : \alpha\in I\}=\{U_\alpha : \alpha<\kappa\}$ and $U_\alpha\neq U_\beta$, provided that $\alpha,\beta\in I$ are different. 

Now, if $I' := \{ \alpha\in I : U_\alpha \neq \emptyset\}$, then clearly $I' \in [I]^{\kappa}$ and hence the hypothesis guarantees the existence of $J\in [I']^{\kappa}$ such that $\{U_\alpha : \alpha\in J\}$ has non-empty intersection (resp., it is centered; it is linked). In these circumstances, a routine argument shows that $J$ and $\{U_\alpha : \alpha\in J\}$ satisfy the desired conditions.
\end{proof}

In what follows we will constantly use the equivalence exposed in Lemma~\ref{lema_equivalencia*}, even without making explicit reference to it.

Let us fix a topological space $X$. The inclusions $\C(X) \subseteq \P(X)\subseteq \WP(X)$ and $\C^*(X) \subseteq \P^*(X)\subseteq \WP^*(X) $ are easily deduced from the definitions. Furthermore, the direct implication of Lemma~\ref{equiv_calibre_indices} indicates that $\C(X) \subseteq \C^*(X)$, $\P(X) \subseteq \P^*(X)$ and $ \WP(X) \subseteq \WP^*(X)$. On the other hand, if $\R$ is the subclass of $\CN$ formed by the regular cardinals, Lemma~\ref{equiv_calibre_indices} ensures that $\C^{*}(X)\cap \R \subseteq \C (X)$, $\P^{*}(X)\cap \R\subseteq \P(X)$ and $\P^{*}(X)\cap \R\subseteq \P(X)$. These lines are the proof of the following result.

\begin{proposition}\label{calibre_vs_calibre*} 
For a topological space $X$ the following inclusions hold.

\begin{enumerate}

\item $\C(X) \subseteq \P(X)\subseteq \WP(X)$ and $\C^*(X) \subseteq \P^*(X)\subseteq \WP^*(X) $.

\item $\C(X) \subseteq \C^*(X)$, $\P(X) \subseteq \P^*(X)$ and $\WP(X) \subseteq \WP^*(X )$.

\item $\C^{*}(X)\cap \R \subseteq \C(X)$, $\P^{*}(X)\cap \R\subseteq \P(X)$ and $\WP^{*}(X)\cap \R\subseteq \WP(X)$.
\end{enumerate}
\end{proposition}

The relations exposed in Proposition~\ref{calibre_vs_calibre*} can be displayed in a more friendly way in the following diagram (the relation $A\to B$ means $A\subseteq B$):

\begin{center}
\begin{tikzcd}

\WP^{*}(X)\cap \R \arrow[r]{}{} & \WP(X) \arrow[r]{}{} & \WP^*(X) \\

\P^{*}(X)\cap \R \arrow[u]{}{} \arrow[r]{}{} & \P(X) \arrow[u]{}{} \arrow[r]{}{} & \P^*(X) \arrow[u]{}{} \\

\C^{*}(X)\cap \R \arrow[u]{}{} \arrow[r]{}{} & \C(X) \arrow[u]{}{} \arrow[r]{}{} & \C^*(X) \arrow[u]{}{} \\

\end{tikzcd}
\end{center}
\vspace{-0.5cm}

The question that we would now like to answer is: Is it true that the inclusions $\C^*(X) \subseteq \C(X)$, $\P^*(X) \subseteq \P(X)$ and $\WP^*(X) \subseteq \WP(X)$ also hold? This question can be answered with the help of the following result.

\begin{proposition}\label{calibre*_topologia} 
If $X$ is a topological space and $\kappa$ is an infinite cardinal such that $o(X)<\kappa$, then $\kappa\in \C^{*}(X)$.
\end{proposition}

\begin{proof} 
If $\kappa$ were not a caliber*, then there would exist a family $\mathcal{U}$ of cardinality $\kappa$ of open subsets of $X$ such that no subcollection of $\mathcal{U}$ of cardinality $\kappa$ has a non-empty intersection, but there can be no such family 
$\mathcal{U}$ since $o(X) < \kappa$.
\end{proof}

The following question naturally arises:

\begin{question}\label{Q1} 
Is it true that if $\kappa\leq o(X)$ and $\kappa\in \C^{*}(X)$, then $\kappa\in\C(X)$?
\end{question}

We will see in Theorems~\ref{thm_1} and \ref{thm_1.2} that, consistently, the answer to Question~\ref{Q1} is in the negative. To do this we first need to prove a couple of auxiliary results.

\begin{lemma}\label{lema_celular_no_WP} 
If $X$ is a topological space, $\kappa$ is an infinite cardinal and $\{U_\alpha : \alpha<\kappa\}$ is a cellular family in $X$, then $\aleph_\kappa \not \in \WP(X)$. Furthermore, if for every $\alpha<\kappa$ it happens that $o(U_\alpha) \geq \aleph_\alpha$, then $\aleph_\kappa \not \in \WP^*(X)$.
\end{lemma}

\begin{proof}
For each $\alpha<\kappa$, consider the set $$I_\alpha := \begin{cases} [0,\aleph_1), & \text{if} \ \alpha=0, \\
[\aleph_{\alpha},\aleph_{\alpha+1}), & \text{if} \ \alpha\neq 0.
   \end{cases}$$ Now, for each ordinal $\gamma<\aleph_\kappa$ define $V_\gamma := U_\alpha$ if and only if $\gamma \in I_\alpha$. It is clear that $\{V_\gamma : \gamma<\aleph_\kappa\}$ is a subset of $\tau_X^{+}$. Moreover, if $J\in [\aleph_\kappa]^{\aleph_\kappa}$, there necessarily exist $\alpha<\beta<\kappa$ such that $I_\alpha\cap J \neq \emptyset\neq I_ \beta\cap J$; consequently, if $\gamma\in I_\alpha\cap J$ and $\delta\in I_\beta\cap J$, then $V_\gamma\cap V_\delta=\emptyset$. Thus, $\aleph_\kappa \not\in \WP(X)$.

 Suppose further that $o(U_\alpha) \geq \aleph_\alpha$, provided that $\alpha<\kappa$. Let us fix a family $\{U(\alpha,\gamma) : \gamma <\aleph_\alpha\} \subseteq \tau_{U_\alpha^+}$, enumerated without repetitions, for every $\alpha<\kappa$. For each ordinal $\gamma<\aleph_\kappa$ define $V_\gamma := U(\alpha,\gamma)$ if and only if $\gamma \in I_\alpha$. Under these circumstances, it is not difficult to confirm that the family $\{V_\gamma : \gamma<\aleph_\kappa\}$ satisfies that $V_\gamma \neq V_\delta$, given that $\gamma<\delta< \aleph_\kappa$. Besides, if $J\in [\aleph_\kappa]^{\aleph_\kappa}$, then there exist $\alpha<\beta<\kappa$ such that $I_\alpha\cap J \neq \emptyset\neq I_ \beta\cap J$. Hence, if $\gamma\in I_\alpha\cap J$ and $\delta\in I_\beta\cap J$, then $V_\gamma\cap V_\delta=\emptyset$. Therefore, $\aleph_\kappa \not\in \WP^*(X)$.
\end{proof}  

Although the notions of caliber and precalibers do not always coincide with their $*$ versions, they do have some similarities. For example, our following result is the natural generalization of Lemma~\ref{precalibre_debil_celularidad} for weak precaliber*.

\begin{lemma}\label{precalibre*_celularidad} 
If $X$ is a topological space, $\kappa$ is a weak precaliber* for $X$ and $\mathcal{U}$ is a cellular family in $X$, then $|\mathcal{U}|<\kappa$.
\end{lemma}

\begin{proof} Let $\mathcal{U}$ be a subset of $\tau_X^+$ such that $|\mathcal{U}|\geq \kappa$. If we take $\mathcal{V}$ in $[\mathcal{U}]^{\kappa}$, we can use that $X$ has weak precaliber* $\kappa$ to find a linked family $\mathcal{W} \in [\mathcal{V}]^{\kappa}$. In particular, $\mathcal{U}$ cannot be a cellular family.
\end{proof}

A natural question suggested by  Lemma~\ref{precalibre*_celularidad} is whether the converse of this result is also valid; that is, is it true that if $X$ is a topological space, 
$\kappa$ is a cardinal number and any cellular family in $X$ has cardinality less than 
$\kappa$, then $\kappa$ is a weak precaliber* for $X$? To answer this question in the affirmative in certain particular cases, we first need to carry out a brief combinatory interlude.

If $\kappa$, $\lambda$, $\mu$ and $\nu$ are cardinal numbers, then $\kappa \to (\lambda)_\nu^\mu$ will mean, as usual, that for any set $X$ with $|X|=\kappa$ and any function $f: [X]^{\mu} \to \nu$, there exists a set $Y\in[X]^{\lambda}$ such that $\restr{f}{[Y]^{\mu}}$ is constant.

A well known combinatory result is Ramsey's Theorem (see \cite[Theorem~10.2, p.~66]{ehmr1984}): 

\begin{theorem}\label{thm_Ramsey} If $r$ and $k$ are positive integers, then $\omega \to \left(\omega\right)_{k}^{r}$. 
\end{theorem}

A cardinal $\kappa>\omega$ is {\it weakly compact} if it satisfies the relation $\kappa \to \left(\kappa\right)_{2}^{2}$. It is in this class of cardinals that we can give a converse of Lemma~\ref{precalibre*_celularidad}. The authors thank Jorge Antonio Cruz Chapital for suggesting the proof that we will present below.

\begin{proposition}\label{equiv_precalibre*_celularidad} 
If $X$ is a topological space and $\kappa=\omega$ or $\kappa$ is weakly compact, then the following statements are equivalent.

\begin{enumerate}
\item $\kappa$ is a weak precaliber* for $X$.
\item If $\mathcal{U}$ is a cellular family in $X$, then $|\mathcal{U}|<\kappa$.
\end{enumerate}
\end{proposition}

\begin{proof} 
The implication $(1)\to (2)$ follows from  Lemma~\ref{precalibre*_celularidad}. To verify the converse suppose that $\kappa$ is not a weak precaliber* for $X$, and let $\mathcal{U} \in [\tau_X^+]^{\kappa}$ such that for any $\mathcal{V}\in [\mathcal{U}]^{\kappa}$, $\mathcal{V}$ is not linked. Let us define a function $f: [\mathcal{U}]^{2} \to 2$ by the rule $f(\{U,V\}) = 0$ if and only if $U\cap V = \emptyset$. Under these circumstances, the relation $\kappa \to \left(\kappa\right)_{2}^{2}$ implies the existence of $\mathcal{V} \in [\mathcal{U}]^{\kappa}$ such that $\restr{f}{[\mathcal{V}]^{2}}$ is constant. Finally, since $\mathcal{V}$ is not linked, $f$ cannot be the constant $1$ in the set $[\mathcal{V}]^{2}$, and therefore we deduce that $\mathcal{ V}$ is a cellular family in $X$ of cardinality $\kappa$.
\end{proof}

\begin{corollary} 
If $X$ is an infinite topological space and $\omega$ is a weak precaliber* for $X$, then $X$ does not have the Hausdorff property.
\end{corollary}

With this background we are better positioned to give a consistent negative answer to Question~\ref{Q1}.

\begin{theorem}\label{thm_1} 
If $\kappa$ is an infinite cardinal with $2^{\kappa}<\aleph_\kappa$ and $\cf(\kappa)>\omega$, then there exists a topological space $X$ such that $X$ is locally compact and $T_4$, $o(X)=2^{\aleph_\kappa}=|X|$, $\aleph_\kappa \in \C^{*}(X)$, $\aleph_\kappa \not\in \WP(X)$ and $\kappa\not \in \WP^*(X)$.
\end{theorem}

\begin{proof}
For each $\xi<\kappa$, let us denote by $Y_\xi$ the Cantor cube of weight $\kappa$, $D(2)^{\kappa}$, and take $Y:=\bigoplus\{ Y_\xi : \xi<\kappa\}$. Also, let $Z$ be the Cantor cube of weight $\aleph_\kappa$, $D(2)^{\aleph_\kappa}$, and define $X:= Y\oplus Z$. Clearly, $X$ is locally compact, $T_4$ and $o(X)=2^{\aleph_\kappa}=|X|$ (see Lemma~\ref{lema_cubo_topologia}). Additionally, by virtue of Lemma~\ref{lema_celular_no_WP}, $\aleph_\kappa \not\in \WP(Y)$ and, since $Y\in \tau_X^+$, we deduce that $\aleph_\kappa \not \in \WP(X)$. On the other hand, since $\{Y_\xi : \xi<\kappa\}$ is a cellular family in $X$, Lemma~\ref{precalibre*_celularidad} implies that $\kappa\not \in \WP ^*(X)$.

To verify that $\aleph_\kappa \in \C^{*}(X)$, suppose that $\mathcal{U}$ is an element of $[\tau_X^+]^{\aleph_\kappa}$ and let $\{U_\alpha : \alpha<\aleph_\kappa\}$ be an enumeration without repetitions of $\mathcal{U}$. For each $\alpha<\aleph_\kappa$, there exist $V_\alpha \in \tau_Y$ and $W_\alpha \in \tau_Z$ such that $U_\alpha = V_\alpha \cup W_\alpha$. We assert that $\left|\{W_\alpha : \alpha<\aleph_\kappa\}\right|=\aleph_\kappa$. Otherwise, given that the function $\mathcal{U} \to \{V_\alpha : \alpha<\aleph_\kappa\} \times \{W_\alpha : \alpha<\aleph_\kappa\} $ determined by $U_\alpha \mapsto (V_\alpha,W_\alpha)$ is injective, we obtain the following contradiction: \[\aleph_\kappa= |\mathcal{U}| \leq \left|\{V_\alpha : \alpha<\aleph_\kappa\}\right|\cdot \left|\{W_\alpha : \alpha<\aleph_\kappa\}\right| < 2^{\kappa} \cdot \aleph_\kappa = \aleph_\kappa.\]

Therefore, there exists $I\in [\aleph_\kappa]^{\aleph_\kappa}$ such that $\emptyset\not\in \{W_\alpha : \alpha \in I\}$ and  $W_\alpha\neq W_\beta$ for any distinct $\alpha,\beta\in I$. Under these circumstances, since $Z$ has caliber $\aleph_\kappa$ (see Lemma~\ref{lema_cubo_calibre}), there exists $J\in [I]^{\aleph_\kappa}$ with $\bigcap\{ W_\alpha : \alpha\in J\}\neq\emptyset$. Consequently, $\{U_\alpha : \alpha\in J\} \in [\mathcal{U}]^{\aleph_\kappa}$ and $\bigcap\{U_\alpha : \alpha\in J\}\neq\emptyset$. Thus, $\aleph_\kappa \in \C^{*}(X)$.
\end{proof}

In Theorem~\ref{thm_1} it might be tempting to consider the softer constraint $2^{\kappa}\leq \aleph_\kappa$, however, since $(2^\kappa)^{\kappa} = 2 ^{\kappa}$ and $(\aleph_\kappa)^{\kappa} > \aleph_\kappa$ (see \cite[Lemma~10.40, p.~34]{kunen1980}), $2^\kappa\neq \aleph_\kappa$ is a result in \textsf{ZFC}.

Moreover, since Theorem~\ref{thm_1} only deals with cardinals with uncountable cofinality, it remains to find out if it is possible to achieve a similar result for uncountable cardinals with countable cofinality. We will settle this debt later in Theorem~\ref{thm_1_p2}.

Since for any $U\in \tau_X$ the inclusion $\tau_U^+ \subseteq \tau_X^+$ is satisfied, the following result is immediate.

\begin{proposition}\label{preservacion*_abiertos} 
Let $X$ be a topological space and $U\in \tau_X$. If $\kappa$ is a caliber* (resp., precaliber*; weak precaliber*) for $X$, then $\kappa$ is a caliber* (resp., precaliber*; weak precaliber*) for $U$.
\end{proposition}

It is easy to verify that the converse implication in the previous result is not true. Indeed, if 
$\kappa$ is an infinite cardinal with $\cf(\kappa)>\omega$, and we denote by $X$ the free sum $\mathbb{R} \oplus D(\kappa)$, then $\mathbb{R} \in \tau_X$ and $\kappa$ is a caliber for $\mathbb{R}$ (see Lemma~\ref{lema_densidad}). However, since $D(\kappa)$ admits a cellular family of cardinality $\kappa$, then Lemma~\ref{precalibre*_celularidad} guarantees that $D(\kappa)$ has no weak precaliber* $\kappa $. Proposition~\ref{preservacion*_abiertos} ensures that $\kappa$ is not a weak precaliber* for $X$.

Recall that a Hausdorff space is {\it $H$-closed} if it is closed in any of its Hausdorff extensions; for example, any compact Hausdorff space is $H$-closed. Naturally, a space is {\it locally 
$H$-closed} if any point admits an $H$-closed neighborhood. A classical Hausdorff extensions result guarantees that if $X$ is a locally $H$-closed Hausdorff space, then 
$X$ is open in any of its Hausdorff extensions (see 
\cite[Proposition~(b), p .~543]{porwoo1988}). 

\begin{corollary}\label{cor_calibres_extensiones} 
If $X$ is a locally $H$-closed Hausdorff space and $Y$ is a $T_2$ extension of $X$, then $\C(X)= \C(Y)$, $\P(X )= \P(Y)$, $\WP(X)= \WP(Y)$, $\C^*(Y)\subseteq \C^*(X)$, 
$\P^*(Y) \subseteq \P^*(X)$ and $\WP^*(Y)\subseteq \WP^*(X)$.
\end{corollary}

\begin{proof} The three equalities follow from Propositions~\ref{calibre_subconjunto_denso} and \ref{calibre_subconjunto_abierto}, while the three inclusions are a consequence of Proposition~\ref{preservacion*_abiertos}.
\end{proof}

Compare the following result with Proposition~\ref{calibre_subconjunto_denso}(2).

\begin{proposition}\label{calibre*_subconjunto_denso} 
If $X$ is a topological space, $D$ is a dense subset of $X$, $\kappa$ is an infinite cardinal, and $\kappa$ is a precaliber* (resp., weak precaliber*) for $X$, then $\kappa$ is a precaliber* (resp., weak precaliber*) for $D$.
\end{proposition}

\begin{proof} 
Let $\{V_\alpha : \alpha<\kappa\}$ be a subset of $\tau_D^{+}$ enumerated without repetitions, and for each $\alpha<\kappa$ set $U_\alpha \in \tau_X^+$ such that $U_\alpha\cap D = V_\alpha$. Then, since the assignment $\alpha\mapsto U_\alpha$ is injective and $X$ has precaliber* (resp., weak precaliber*) $\kappa$, there exists $J\in [\kappa]^{\kappa} $ such that $\{U_\alpha : \alpha\in J\}$ is centered (resp., linked). Finally, the density of $D$ guarantees that $\{V_\alpha : \alpha\in J\}$ is centered (resp., linked).
\end{proof}

The corresponding result for calibers* is invalid. Indeed, Lemma~\ref{lema_cubo_calibre} guarantees that if $\kappa$ has uncountable cofinality, then $D(2)^{\kappa}$ has caliber $\kappa$. On the other hand, a routine argument shows that $D := \{x\in D(2)^{\kappa} : 
|\{\alpha<\kappa : x(\alpha)\neq 0\}|\leq \omega\}$ is a dense subset of $D(2)^{\kappa}$. Furthermore, if for each $\alpha<\kappa$ we denote by $\pi_\alpha : D(2)^{\kappa} \to D(2)$ the $\alpha^{\text{th}}$ canonical projection, then $\{\pi_\alpha^ {-1}\{1\}\cap D : \alpha<\kappa\}$ is a family of $\kappa$ open subsets such that for any $J\in [\kappa]^{\kappa}$, $\bigcap\{\pi_\alpha^{-1}\{1\}\cap D : \alpha\in J\}=\emptyset$. Consequently, $\kappa$ is not a caliber* for $D$.

\begin{question}\label{Q2} 
Is it true that $\kappa$ is a caliber* (resp., precaliber*; weak precaliber*) for a topological space $X$, if $\kappa$ is a caliber* (resp., precaliber*; weak precaliber*) for a dense subset $D$ of $X$?
\end{question}

We will show in Corollaries~\ref{cor1_thm_extension} and \ref{cor2_thm_extension} that the answer to Question~\ref{Q2} is in the negative. However, we will first discuss some simple cases in which the answer is in the affirmative. For instance, the following result is a consequence of Proposition~\ref{calibre_subconjunto_denso} and Lemma~\ref{equiv_calibre_indices}.

\begin{proposition} If $X$ is a topological space, $D$ a dense subset of $X$, and $\kappa$ a regular cardinal that is a caliber* (resp., precaliber*; weak precaliber*) for $D$, then $\kappa$ is a caliber* (resp., precaliber*; weak precaliber*) for $X$.
\end{proposition}

Corollary~\ref{cor_loc_comp_caliber*} below shows that in a certain class of extensions the answer to Question~\ref{Q2} is in the affirmative. The following result is a consequence of Proposition~\ref{preservacion*_union_regular} and Theorem~\ref{preservacion*_union} that we will prove later.

\begin{proposition}\label{preservacion*_residuo_calibre} 
Let $X$ be a topological space and $Y$ an extension of $X$. If $\kappa$ is a caliber* (resp., precaliber*; weak precaliber*) for $X$ and for the remainder $Y\setminus X$, then $\kappa$ is a caliber* (resp., precaliber*; weak precaliber*) for $Y$.
\end{proposition}

Moreover, it is possible to give a kind of reciprocal result for the above by virtue of Proposition~\ref{preservacion*_abiertos}.

\begin{proposition}\label{preservacion*_residuo_cerrado} 
Let $X$ be a topological space and $Y$ an extension of $X$. If $\kappa$ is a caliber* (resp., precaliber*; weak precaliber*) for $Y$, and $Y\setminus X$ is closed in $Y$, then $\kappa$ is a caliber* (resp., precaliber*; weak precaliber*) for $X$. 
\end{proposition}

Propositions~\ref{preservacion*_residuo_calibre} and \ref{preservacion*_residuo_cerrado} imply the following corollary related with the Alexandroff one-point compactification $\alpha X$ of a topological
space $X$.

\begin{corollary}\label{cor_loc_comp_caliber*} If $X$ is locally compact, non-compact and Hausdorff, then $\C^{*}(X) = \C^{*}(\alpha X)$, $\P^{*}(X) = \P^{*}(\alpha X)$ and $\WP^{*}(X) = \WP^{*}(\alpha X)$.
\end{corollary}

With the previous result available it is possible to strengthen Theorem~\ref{thm_1} as follows:

\begin{theorem}\label{thm_1.2} If $\kappa$ is an infinite cardinal with $2^{\kappa}<\aleph_\kappa$ and $\cf(\kappa)>\omega$, then there exists a compact Hausdorff space $X$ such that $o(X)=2^{\aleph_\kappa}=|X|$, $\aleph_\kappa \in \C^{*}(X)$, $\aleph_\kappa \not\in \WP(X)$ and $\kappa\not \in \WP^*(X)$.

\end{theorem}

\begin{proof} Let $X$ be the one-point compactification of the space whose existence is guaranteed by Theorem~\ref{thm_1}. Clearly, $X$ is a compact Hausdorff space that satisfies the relations $o(X)=2^{\aleph_\kappa}=|X|$. Furthermore, Proposition~\ref{calibre_subconjunto_denso}(1) confirms that $\aleph_\kappa \not\in \WP(X)$, while Corollary ~\ref{cor_loc_comp_caliber*} ensures that $\aleph_\kappa \in \C^{*}(X)$ and $\kappa\not \in \WP^*(X)$.
\end{proof}

\begin{question} Is it possible to construct in \textsf{ZFC} an example of an infinite Hausdorff space $X$ and an infinite cardinal $\kappa$ such that $\kappa\in \C^*(X)$ and $\kappa \not\in \C(X)$?

\end{question}

In order to give a negative answer to Question~\ref{Q2}, let us first prove the following theorem.

\begin{theorem}\label{thm_extension} 
If $\kappa$ is a singular cardinal and $X$ is a topological space such that $\cf(\kappa)$ is not a caliber* (resp., precaliber*; weak precaliber*) for $X$, then there exists an extension $Y$ of $X$ such that $\kappa$ is not a caliber* (resp., precaliber*; weak precaliber*) for $Y$.
\end{theorem}

\begin{proof} 
Let $\{\kappa_\xi : \xi<\cf(\kappa)\}$ be a strictly increasing sequence of cardinal numbers with $\kappa_0 = 0$ and $\sup\{\kappa_\xi : \xi<\cf(\kappa)\} =\kappa$. Also, for each $\xi<\cf(\kappa)$, let us define $I_\xi := [\kappa_\xi, \kappa_{\xi+1})$. Take a collection $\{y_\alpha : \alpha<\kappa\}$ such that $X \cap \{y_\alpha : \alpha<\kappa\} = \emptyset$ and $y_\alpha \neq y_\beta$, provided that $\alpha<\beta<\kappa$. Let us use the fact that $\cf(\kappa)$ is not a caliber* (resp., precaliber*; weak precaliber*) for $X$ to fix a family $\{V_\xi :  \xi<\cf(\kappa) \}\subseteq \tau_X^+$ enumerated without repetitions and such that for any $I\in [\cf(\kappa)]^{\cf(\kappa)}$, the collection $\{V_\xi : \xi\in I\}$ has an empty intersection (resp., it is not centered; it is not linked). Additionally, for each $\alpha<\kappa$ define $U_\alpha := V_\xi \cup \{y_\alpha\}$ if and only if $\alpha \in I_\xi$. Finally, let $Y := X\cup \{y_\alpha : \alpha<\kappa\}$ and consider the topology on $Y$ generated by the base $\tau_X \cup \{U_\alpha : \alpha <\kappa\}$.

Clearly, $X$ is a dense subspace of $Y$. Note also that if $\alpha<\beta<\kappa$ and 
$\xi,\eta<\cf(\kappa)$ are such that $\alpha \in I_\xi$ and $\beta \in I_\eta$, then $U_\alpha\cap U_\beta = V_\xi \cap V_\eta$. Now, let $J\in [\kappa]^{\kappa}$ and consider the collection $I:= \{\xi<\cf(\kappa) : I_\xi \cap J \neq \emptyset\} $. First note that since 
$|J|=\kappa$, necessarily $|I|=\cf(\kappa)$; in particular, $\{V_\xi : \xi\in I\}$ has an empty intersection (resp., is not centered; is not linked). Under these circumstances, $\{U_\alpha : \alpha\in J\}$ has an empty intersection (resp., is not centered; is not linked). Thus, $\kappa$ is not a caliber* (resp., is not a precaliber*; is not a weak precaliber*) for $Y$. 
\end{proof}

\begin{corollary}\label{cor1_thm_extension} 
If $\kappa$ is a regular uncountable cardinal with $2^{\kappa}<\aleph_\kappa$, then there exist a topological space $X$ and an extension $Y$ of $X$ such that $\aleph_\kappa \in \C^{*}(X)$ and $\aleph_\kappa \not\in \WP^{*}(Y)$.
\end{corollary}

\begin{proof} Let $X$ be as in Theorem~\ref{thm_1}. Since $\aleph_\kappa$ is a singular cardinal and $\cf\left(\aleph_\kappa\right) = \kappa$ is not a weak precaliber* for $X$, Theorem~\ref{thm_extension} produces an extension $Y$ of $X$ such that $\aleph_\kappa$ is not a weak precaliber* for $Y$. Thus, $\aleph_\kappa \in \C^{*}(X)$ and $\aleph_\kappa \not\in \WP^{*}(Y)$.
\end{proof}

In general, it is not a simple task to determine precisely what the chain conditions for topological spaces are. In what follows we will try to exemplify this fact by means of some well known spaces.

Recall that an infinite family $\mathcal{A}$ formed by infinite subsets of $\omega$ is {\it almost disjoint} if $|A\cap B|<\omega$, provided that $A,B\in\mathcal{A}$ are different. A classical result guarantees the existence of an almost disjoint family of cardinality 
$\mathfrak{c}$ (see, for example, \cite[Theorem~1.3, p.~48]{kunen1980}). In particular, for any cardinal $\kappa$ with $\omega\leq \kappa\leq\mathfrak{c}$, there exists an almost disjoint family of cardinality $\kappa$ on $\omega$.

\begin{lemma}\label{lema_familia_casi_ajena} 
If $\kappa\geq\omega$ is a cardinal number with $\kappa<\mathfrak{c}$ and $\cf(\kappa)=\omega$, then there exists an almost disjoint family $\mathcal{A}$ on $\omega$, such that $|\mathcal{A}|=\kappa$ and for any $\mathcal{B} \in [\mathcal{A}]^{\kappa}$, there exist $A,B,C\in\mathcal{B}$ with $A\cap B \cap C = \emptyset$.
\end{lemma}

\begin{proof} 
Let $\{\kappa_n : n<\omega\}$ be a strictly increasing sequence of cardinal numbers with $\kappa_0 = 0$ and $\sup\{\kappa_n : n<\omega\} =\kappa$. Also, for each $n<\omega$, let us define $I_n := [\kappa_n, \kappa_{n+1})$. Furthermore, let us use the relation $\kappa<\mathfrak{c}$ to fix an almost disjoint family $\{B_\alpha : \alpha<\kappa\}$ enumerated without repetitions. The last step is to define for any $\alpha<\kappa$, $A_\alpha := B_\alpha \setminus n$ if and only if $\alpha \in I_n$.

\medskip

\noindent {\bf Claim.} The collection $\{A_\alpha : \alpha<\kappa\}$ is a faithfully indexed almost disjoint family such that, for any $J\in [\kappa]^{\kappa}$, there exist $ \alpha,\beta,\gamma\in J$ with $A_\alpha\cap A_\beta\cap A_\gamma=\emptyset$.

\medskip

Let us first notice that, if $\alpha<\beta<\kappa$, and $n,m\in \omega$ are such that $\alpha\in I_n$ and $\beta\in I_m$, then the relation $A_ \alpha = A_\beta$ implies that $B_\alpha \setminus n = B_\beta \setminus m$. Consequently, $B_\alpha$ is contained in the union $n\cup B_\beta$ and hence $B_\alpha = (B_\alpha\cap n) \cup (B_\alpha\cap B_\beta )$, which is absurd since $B_\alpha$ is infinite, while this last union is finite. For this reason, $A_\alpha \neq A_\beta$. Also, since $|A_\alpha \cap A_\beta| \leq |B_\alpha \cap B_\beta| <\omega$, and $\{A_\gamma : \gamma<\kappa\}$ is clearly formed by infinite subsets of $\omega$, it follows that $\{A_\gamma : \gamma< \kappa\}$ is an almost disjoint family enumerated without repetitions.

On the other hand, if $J\in [\kappa]^{\kappa}$ and $\alpha,\beta \in J$ are different, there exists $m<\omega$ with $A_\alpha \cap A_\beta \subseteq m$. Thus, if we take $n>m$ with $I_n \cap J\neq \emptyset$ and $\gamma\in I_n \cap J$, then the relations $A_\alpha \cap A_\beta \cap A_\gamma \subseteq m \cap (B_\gamma \setminus n) = \emptyset$ hold. Consequently, $A_\alpha \cap A_\beta \cap A_\gamma = \emptyset$.

Finally, the Claim guarantees that $\mathcal{A} := \{A_\alpha : \alpha<\kappa\}$ is the family required in the statement of the present lemma.
\end{proof}

Recall that the Stone-\v{C}ech compactification of $\omega$, $\beta\omega$, is the space of ultrafilters on $\omega$ equipped with the topology generated by the base 
$\{A^* : A\subseteq \omega\}$ where for each $A\subseteq \omega$, $A^* := \{\mathcal{U} \in \beta\omega : A\in\mathcal{U}\}$.

\begin{remark}\label{obs_beta_omega} 
If $\mathcal{A}$ is an almost disjoint family on $\omega$ and $\mathcal{U} \in\beta\omega\setminus\omega$, then $|\mathcal{A} \cap \mathcal{ U}|\leq 1$. Indeed, if $A,B\in\mathcal{A}$ are different, then $|A\cap B|<\omega$ and thus, since $\mathcal{U}$ is free, we deduce that $A\cap B \not \in \mathcal{U}$; consequently, $\{A,B\} \not \subseteq \mathcal{U}$.
\end{remark}

\begin{lemma}\label{lema_beta_omega} 
If $\kappa\geq\omega$ is a cardinal number with $\kappa<\mathfrak{c}$ and $\cf(\kappa)=\omega$, then $\kappa\not\in \P^{* }(\beta\omega)$.
\end{lemma}

\begin{proof} Let $\mathcal{A}$ be an almost disjoint family with the characteristics of Lemma~\ref{lema_familia_casi_ajena}, and consider the family $\mathbb{A}:=\{A^* : A\in\mathcal{A} \}$. Since each element of $\mathcal{A}$ is contained in a free ultrafilter on $\omega$, Remark~\ref{obs_beta_omega} implies that $\mathbb{A}\in [\tau_{\beta\omega}^+]^{\kappa}$. Now let $\mathbb{B}\in [\mathbb{A}]^{\kappa}$, $\mathcal{B} \in [\mathcal{A}]^{\kappa}$ where $\mathbb {B} = \{B^* : B\in \mathcal{B}\}$, and $A,B,C\in\mathcal{B}$ such that $A\cap B \cap C = \emptyset $. Thus, there is no $\mathcal{U} \in \beta\omega$ such that $\{A,B,C\} \subseteq \mathcal{U}$; that is, $A^* \cap B^* \cap C^*=\emptyset$. Hence, $\mathbb{A}$ is an element of $[\tau_{\beta\omega}^+]^{\kappa}$ and, for any $\mathbb{B}\in [\mathbb{A} ]^{\kappa}$, $\mathbb{B}$ is not a centered family. Consequently, $\kappa$ is not a precaliber* for $\beta\omega$.
\end{proof}

\begin{theorem} 
The following statements are true.

\begin{enumerate}
\item $\C(\beta\omega) = \P(\beta\omega) = \WP(\beta\omega) = \{\kappa\in\CN : \cf(\kappa)>\omega\} $.
\item $\{\kappa\in\CN : \cf(\kappa)>\omega\}\cup \{\kappa\in\CN : \kappa\geq 2^{\mathfrak{c}}\} \subseteq \C^{*}(\beta\omega)$.
\item $\P^{*}(\beta\omega) \subseteq \{\kappa\in\CN : \cf(\kappa)>\omega\}\cup \{\kappa\in\CN : \cf (\kappa)=\omega \ \wedge \ \kappa>\mathfrak{c}\}$.
\end{enumerate}

In particular, under the hypothesis $2^{\mathfrak{c}} = \mathfrak{c}^+$,

\begin{enumerate}
\item[(4)] $\C^{*}(\beta\omega) = \P^{*}(\beta\omega) = \{\kappa\in\CN : \kappa\geq 
2^{\mathfrak{c}}\}\cup \{\kappa\in\CN : \kappa< 2^{\mathfrak{c}}\ \wedge \ \cf(\kappa)>\omega\}$.
\end{enumerate}
\end{theorem}

\begin{proof} 
First, since $\beta\omega$ is a separable Hausdorff space, a combination of Lemmas~\ref{lema_densidad}, \ref{precalibre_debil_celularidad}, and \ref{calibre_subconjunto_denso} implies that $\C(\beta\omega) = \P(\beta\omega) = \WP(\beta\omega) = \{\kappa\in\CN : \cf(\kappa)>\omega\}$.

On the other hand, by virtue of item (1) of this theorem and  Proposition~\ref{calibre*_topologia}, the inclusion $\{\kappa\in\CN : \cf(\kappa)>\omega\}\cup \{\kappa\in\CN : \kappa\geq 2^{\mathfrak{c}}\} \subseteq \C^{*}(\beta\omega)$ is clear. Moreover, if $\kappa$ is an infinite cardinal with $\kappa<\mathfrak{c}$ and $\cf(\kappa)=\omega$, then  Lemma~\ref{lema_beta_omega} ensures that $\kappa\not\in\P^{*}(\beta\omega)$. Consequently, $\P^{*}(\beta\omega) \subseteq \{\kappa\in\CN : \cf(\kappa)>\omega\}\cup \{\kappa\in\CN : \cf(\kappa)=\omega \ \wedge \ \kappa>\mathfrak{c}\}$.
\end{proof}

\begin{question}\label{Q_beta_omega1} What can we say about the collection
$\WP^*(\beta\omega)$?
\end{question}

In particular, certain chain conditions depend entirely on how the cardinal numbers are related to each other.

\begin{corollary} 
The following statements are true.

\begin{enumerate}
\item $\aleph_\omega<\mathfrak{c}$ implies that $\aleph_\omega \not \in \P^{*}(\beta\omega)$.
\item $\aleph_\omega>2^{\mathfrak{c}}$ implies that $\aleph_\omega \in \C^{*}(\beta\omega)$.
\end{enumerate}
\end{corollary}

\begin{question}\label{Q_beta_omega2} 
Is it true that if $\kappa$ is a cardinal with $\mathfrak{c}<\kappa<2^{\mathfrak{c}}$ and $\cf(\kappa)=\omega$, then $\kappa$ is a caliber* for $\beta\omega$?
\end{question}

\begin{example}\label{ejem_discreto} 
If $\kappa$ is an infinite cardinal, then:

\begin{enumerate}
\item $\C\left(D(\kappa)\right)=\P\left(D(\kappa)\right)=\WP\left(D(\kappa)\right)=\{\lambda\in\CN : \cf(\lambda)>\kappa\}$; and
\item $\{\lambda\in\CN : \lambda >2^{\kappa}\} \subseteq \C^{*}\left(D(\kappa)\right) \subseteq \P^{*} \left(D(\kappa)\right) \subseteq \WP^{*}\left(D(\kappa)\right) \subseteq \{\lambda\in\CN : \lambda >\kappa\}$.
\end{enumerate}

Also, if $2^{\kappa}=\kappa^+$,

\begin{enumerate}
\item[(3)] $\C^{*}\left(D(\kappa)\right)=\P^{*}\left(D(\kappa)\right)=
\WP^{*}\left(D(\kappa)\right) =\{\lambda\in\CN : \lambda >\kappa\}$.
\end{enumerate}
\end{example}

\begin{proof} 
First, since $d\left(D(\kappa)\right)=\kappa$,  Lemma~\ref{lema_densidad} guarantees that the inclusion 
$\{\lambda\in\CN : \cf(\lambda)> \kappa\}\subseteq \C\left(D(\kappa)\right)$ holds. On the other hand, if $\cf(\lambda)\leq \kappa$, then the collection $\{\{\alpha\} : \alpha<\cf(\lambda)\}$ is a cellular family in $D (\kappa)$ of cardinality $\cf(\lambda)$. Thus, Lemmas~\ref{lema_cofinalidad} and \ref{precalibre_debil_celularidad} imply that $\lambda$ is not a weak precaliber for $D(\kappa)$ and therefore the inclusion $\WP\left(D(\kappa)\right) \subseteq \{\lambda\in\CN : \cf(\lambda)>\kappa\}$ is satisfied. Finally, as the relations $\C\left(D(\kappa)\right) \subseteq \P\left(D(\kappa)\right)\subseteq \WP\left(D(\kappa)\right)$ always hold, (1) is true.

It follows from Lemma~\ref{precalibre*_celularidad} that cardinals less than or equal to $\kappa$ cannot be weak precalibers* for $D(\kappa)$ and therefore, $\WP^{*}\left(D(\kappa)\right) \subseteq \{\lambda\in\CN : \lambda >\kappa\}$. Furthermore, since $|\tau_{D(\kappa)}|=2^{\kappa}$,  Proposition~\ref{calibre*_topologia} guarantees that 
$\{\lambda\in\CN : \lambda >2 ^{\kappa}\} \subseteq \C^{*}\left(D(\kappa)\right)$, which proves (2).

Finally, to prove (3), note that since $\kappa^+$ is a regular cardinal greater than $\kappa$, (1), together with (2) of Proposition~\ref{calibre_vs_calibre*}, shows that $\kappa^ {+}\in \C^*\left(D(\kappa)\right)$. If we assume that $2^{\kappa}=\kappa^+$, then the relations in item (2) allow us to see that: \begin{align*} \{\lambda\in\CN : \lambda >\kappa\} &= \{\lambda\in\CN : \lambda \geq \kappa^+\} = \{\kappa^+\}\cup \{\lambda\in\CN : \lambda >2^{\kappa}\} \\
&\subseteq \C^{*}\left(D(\kappa)\right) \subseteq \P^{*}\left(D(\kappa)\right) \\
&\subseteq \WP^{*}\left(D(\kappa)\right) \subseteq \{\lambda\in\CN : \lambda >\kappa\};
\end{align*} therefore, item (3) is satisfied.
\end{proof}

For our next example it is convenient to establish a couple of auxiliary results. Recall that a topological space $X$ is \textit{hyperconnected} if any two non-empty open subsets of $X$ intersect. In other words, $X$ is hyperconnected if and only if $\tau_X^+$ is a linked family. It turns out that this condition is equivalent to a property that is seemingly stronger.

\begin{proposition}\label{equiv_hiperconexo} 
A topological space is hyperconnected if and only if $\tau_X^+$ is a centered family.
\end{proposition}

\begin{proof} 
The inverse implication is clear. To verify the direct implication, we proceed by finite induction. Suppose that for an $n<\omega$ it has already been shown that the intersection of any $n$ elements of $\tau_X^+$ is a non-empty set. Now suppose that $\{U_k : k<n+1\}$ is a subset of $\tau_X^+$ of size $n+1$. It then follows that $\bigcap\{U_k : k<n\}$ is an element of $\tau_X^+$ and therefore, $\bigcap\{U_k : k<n+1\}$ is non-empty.
\end{proof}

From the previous proposition, two more conditions equivalent to being hyperconnected can be obtained in terms of precalibers and weak precalibers. For the following result only, we naturally extend the concepts of precaliber, precaliber*, weak precaliber, and weak precaliber* for any finite cardinal. Note that $1$ is always a precaliber of any space $X$.

\begin{proposition}\label{prop_hiperconexo} 
Let us consider the following statements for a topological space $X$.

\begin{enumerate}
\item $X$ is hyperconnected.

\item Any cardinal number $\kappa > 0$ is a precaliber for $X$.

\item There is a finite cardinal $k > 1$ that is a precaliber of $X$.

\item Any cardinal number $\kappa > 0$ is a weak precaliber of $X$.

\item There is a finite cardinal $k > 1$ that is a weak precaliber for $X$.

\item Any cardinal number $\kappa > 0$ is precaliber* for $X$.

\item There is a finite cardinal $k > 1$ that is a precaliber* of $X$.

\item Any cardinal number $\kappa > 0$ is a weak precaliber* of $X$.

\item There is a finite cardinal $k > 1$ that is a weak precaliber* for $X$.
\end{enumerate}

Then the following statements are verified.

\begin{enumerate}
\item[(a)] $(1)-(5)$ are equivalent.
\item[(b)] $(2)\to(6)$, $(3)\to(7)$, $(4)\to(8)$ and $(5)\to(9)$.
\item[(c)] $(6)\to (7)$, $(6)\to (8)$, $(8) \to (9)$ and $(7)\to (9)$.
\item[(d)] (9) implies (1) when $o(X)\geq \omega$.
\end{enumerate}
In particular, if $\tau_X$ is infinite, then $(1)-(9)$ are equivalent.
\end{proposition}

\begin{proof} 
For item (a), given that implications $(2)\to (3)$, $(2)\to (4)$, $(4) \to (5)$ and $(3) \to (5)$ are clear, we only have to argue implications $(1)\to (2)$ and $(5)\to (1)$.

To verify that (1) implies (2), let us fix a non-zero cardinal number $\kappa$. Since $X$ is hyperconnected,  Proposition~\ref{equiv_hiperconexo} implies that $\tau_X^+$ is a centered family. Thus, if $\{U_\alpha : \alpha < \kappa\}$ is a subset of $\tau_X^+$ indexed by $\kappa$, it must be centered. Therefore $\kappa$ is a precalibre of $X$.

To see that (5) implies (1), suppose that $1<k<\omega$ is a weak precaliber of 
$X$. Let $U_0$ and $U_1$ be two non-empty open subsets of $X$. Let $U_i = U_0$ for each $2\leq i<k$. Since $k$ is a precaliber of $X$, then there exists $J \in [k]^{k}$ such that $\{U_i : i \in J\}$ is linked. But the only subset of $k$ of cardinality $k$ is itself. Therefore $\{U_i : i <k\}$ is linked. In particular $U_0 \cap U_1 \not= \emptyset$. That is, every two non-empty open subsets of $X$ have a non-empty intersection, i.e., $X$ is hyperconnected.

Now, item (b) is a consequence of  Proposition~\ref{calibre_vs_calibre*}.(2), while item (c) is clear. Lastly, for item (d), if $k>1$ is a finite cardinal such that $k$ is a weak precaliber* for $X$, and $U,V\in \tau_X^+$ are different, then there exists $\mathcal{U} \subseteq \tau_X^{+}$ with $|\mathcal{U}\cup \{U,V\}|=k$. Then, since the elements of 
$\mathcal{U}\cup \{U,V\}$ have pairwise non-empty intersection, $U\cap V\neq \emptyset$, which shows that $X$ is hyperconnected.
\end{proof}

It remains to comment on whether the restriction that we had to add in item (d) of the previous result is essential. Clearly, it is satisfied that $k$ is a weak precaliber* for $D(2)$ for any $k\geq 4$, but $D(2)$ is not hyperconnected. On the other hand, a routine argument shows that $X$ is a hyperconnected space if and only if $2$ is a weak precaliber* for $X$. Let us now consider the condition that 3 is a weak precaliber* for $X$.

Note first that if $k$ and $l$ are a pair of finite cardinals such that $l \leq k <o(X)$, then if $k$ is a weak precaliber* for $X$, $l$ is also is a weak precaliber* for $X$. Indeed, if 
$\mathcal{U}\in [\tau_X^+]^{l}$ and $\mathcal{V} \subseteq \tau_X^+$ are such that 
$|\mathcal{U}\cup \mathcal{V}|=k$, then, since our hypothesis guarantees that 
$\mathcal{U}\cup \mathcal{V}$ is a linked family, we deduce that $\mathcal{U}$ also is. In particular, if a finite $k$ with $2\leq k<o(X)$ is a weak precaliber* for $X$, then $X$ is a hyperconnected space.

Now let us assume that $3$ is a weak precaliber* for a space $X$. When $3<o(X)$, the previous paragraph implies that $X$ is hyperconnected; otherwise, the condition $o(X)\leq 3$ implies that either $\tau_X = \{\emptyset, X\}$ or there exists $A\in P(X)\setminus\{\emptyset, X\}$ such that $\tau_X = \{\emptyset, X, A\}$. In either case the space $X$ is hyperconnected.

The classical examples of hyperconnected spaces are the infinite spaces equipped with the cofinite topology. What we will do next is calculate all the above chain conditions for this class of spaces. Recall that an infinite space $X$ is called {\it cofinite} if $\tau_X = \{\emptyset\}\cup \{U\subseteq X : \abs{X\setminus U}<\omega\}$. The following observation is essential.

\begin{remark}\label{obs_calibre_cofinito} 
If $X$ is a cofinite space and $\kappa$ is a cardinal number, then $\kappa$ is a caliber for $X$ if and only if for any $\{F_\alpha : \alpha<\kappa\} \subseteq [X]^{<\omega}$ there exists $J\in [\kappa]^{\kappa}$ such that $\bigcup\{F_\alpha : \alpha\in J\} \neq X $.
\end{remark}

\begin{lemma}\label{lema_calibres_cofinitos} If $X$ is the cofinite space of cardinality $\kappa$, then:
\begin{enumerate}
\item $\P(X)=\P^*(X)=\WP(X)=\WP^*(X)=\CN$;
\item $\left\{\lambda \in\CN : \kappa<\lambda\right\} \subseteq \C^{*}(X)$; and
\item $\{\lambda\in\CN : \lambda < \kappa\ \vee\ \cf(\lambda)>\omega\}\cup\{\lambda\in\CN : \cf(\lambda)<\kappa<\lambda\} \subseteq \C(X)$.
\end{enumerate}
\end{lemma}

\begin{proof} First, since $X$ is a hyperconnected topological space, Proposition~\ref{prop_hiperconexo} ensures that $\CN \subseteq \P(X)$. Consequently, item (1) is true. Item (2) is a consequence of the equality $o(X)=\kappa$ and Proposition~\ref{calibre*_topologia}.

To verify item (3) let us start by observing that, since each element of $[X]^{\omega}$ is dense in $X$, it follows that $d(X)=\omega$. Thus, Lemma~\ref{lema_densidad} guarantees  $\{\lambda\in\CN : \cf(\lambda)>\omega\} \subseteq \C(X)$.

Now, suppose $\lambda\in\CN$ is such that $\lambda<\kappa$ and take a subset $\{F_\alpha : \alpha<\lambda\}$ of $[X]^{<\omega}$. By virtue of the inequalities $\left|\bigcup\left\{F_\alpha : \alpha<\lambda\right\}\right| \leq \lambda\cdot\sup\left\{\left|F_\alpha\right| : \alpha<\lambda\right\}\leq \lambda\cdot \omega <\kappa$, we have that 
$\bigcup\left\{F_\alpha : \alpha<\lambda\right\} \neq X$. Thus, Remark~\ref{obs_calibre_cofinito} ensures that $\lambda \in \C(X)$.

To verify the inclusion $\{\lambda\in\CN : \cf(\lambda)<\kappa<\lambda\} \subseteq \C(X)$, let $\lambda$ be a cardinal number with $\cf( \lambda)<\kappa<\lambda$. Fix a strictly increasing sequence $\{\lambda_\xi : \xi<\cf(\lambda)\}$ such that $\sup\{\lambda_\xi : \xi<\cf(\lambda)\} = \lambda$. With  Remark~\ref{obs_calibre_cofinito} in mind, let $\{F_\alpha : \alpha<\lambda\}$ be a subset of $[X]^{<\omega}$ and consider the function 
$f: \lambda \to [X]^{<\omega}$ defined as $f(\alpha) := F_\alpha$.

Let us construct, through transfinite recursion on $\cf(\lambda)$, a sequence of ordinal numbers $\{\alpha_\xi : \xi<\cf(\lambda)\}$ such that for each $\xi<\cf(\lambda)$, 
$|f^{-1}\{F_{\alpha_\xi}\}|\geq \lambda_\xi$. Suppose that for some $\xi<\cf(\lambda)$ we have constructed a sequence $\{\alpha_\eta : \eta<\xi\}$ with the desired characteristics. If there is no $\alpha<\lambda$ with  $|f^{-1}\{F_{\alpha}\}|\geq \lambda_\xi$, then $$\lambda = \left|\bigcup_{\alpha<\lambda} f^{-1}\{F_\alpha\}\right| \leq \left|[X]^{<\omega}\right| \cdot \sup \left\{\left|f^{-1}\{F_\alpha\}\right| : \alpha<\lambda\right\}\leq \kappa\cdot \lambda_\xi<\lambda;$$ which is a contradiction. Consequently, there is $\xi<\cf(\lambda)$ with $|f^{-1}\{F_{\alpha_\xi}\}|\geq \lambda_\xi$ and hence $\{\alpha_\eta : \eta<\xi+1\}$ meets the required properties.

Now define $J := \bigcup\{f^{-1}\{F_{\alpha_\xi}\} : \xi<\cf(\lambda)\}$ and observe that the properties of $\{\alpha_ \xi : \xi<\cf(\lambda)\}$ guarantee the equality $|J| = \lambda$. Finally, since $\bigcup\{F_\alpha : \alpha\in J\} = \bigcup\{F_{\alpha_\xi}: \xi<\cf(\lambda)\}$ and $$\left| \bigcup\{F_{\alpha_\xi}: \xi<\cf(\lambda)\}\right| \leq \cf(\lambda)\cdot \sup\left\{\left|F_{\alpha_\xi}\right| : \xi<\cf(\lambda)\right\} \leq \cf(\lambda) \cdot \omega<\kappa,$$ we deduce that $\bigcup\{F_\alpha : \alpha\in J\} \neq X$.
\end{proof}

The above result allows us to fully classify the chain conditions on cofinite spaces.

\begin{proposition}\label{prop_cofinito_omega} If $X$ is the cofinite space of cardinality $\omega$, then:

\begin{enumerate}
\item $\P(X)=\P^*(X)=\WP(X)=\WP^*(X)=\CN$;
\item $\C^{*}(X)=\{\lambda \in\CN : \lambda>\omega\}$; and
\item $\C(X)=\{\lambda \in\CN : \cf(\lambda)>\omega\}$.
\end{enumerate}

\end{proposition}

\begin{proof} 
By virtue of Lemma~\ref{lema_calibres_cofinitos}, it is only necessary to verify that $\omega\not\in \C^{*}(X)$ and $\C(X)\subseteq\{\lambda \in\CN : \cf(\lambda)>\omega\}$. If $X=\{x_n : n<\omega\}$ and for each $n<\omega$ we define $F_n:= \{x_k : k\leq n\}$, then $\{F_n : n<\omega\}\subseteq [X]^{<\omega}$ and for each $J\in [\omega]^{\omega}$, $\bigcup\{F_n : n\in J\}=X$. Consequently, Remark~\ref{obs_calibre_cofinito} ensures that $\omega\not\in \C(X)$ and, therefore, $\omega\not\in \C^{*}(X)$ (see  Lemma~\ref{equiv_calibre_indices}). Finally, the relation $\omega\not\in \C(X)$ and Lemma~\ref{lema_cofinalidad} imply that $\C(X)\subseteq \{\lambda\geq\omega : \cf(\lambda)>\omega\}$.
\end{proof}

\begin{proposition}\label{prop_cofinito_>omega} 
If $X$ is the cofinite space of cardinality $\kappa>\omega$, then $\C(X)=\C^{*}(X)=\P(X)=\P^*(X)= \WP(X)=\WP^*(X)=\CN$.
\end{proposition}

\begin{proof} 
Clearly, it is sufficient to verify that $\CN \subseteq \C(X)$. With this objective in mind, we will first argue that $\kappa\in \C(X)$.

If $\cf(\kappa)>\omega$, the relation $\kappa\in \C(X)$ is a consequence of Lemma~\ref{lema_calibres_cofinitos}(3). When $\cf(\kappa)=\omega$, we take a subset $\{F_\alpha : \alpha<\kappa\}$ of $[X]^{<\omega}$ and fix a subset $Y$ of $X$ with $|Y|=\aleph_1$. Consider $Y$ as a subspace of $X$; then $\{F_\alpha \cap Y : \alpha<\kappa\}$ is a subset of $[Y]^{<\omega}$. Since, as $\kappa\in \C(Y)$ (see Lemma~\ref{lema_calibres_cofinitos}(3)), we obtain the existence of $J\in [\kappa]^{\kappa}$ in such a way that $\bigcup\{F_\alpha \cap Y : \alpha\in J\} \neq Y$. Consequently, $\bigcup\{F_\alpha : \alpha\in J\} \neq X$, and therefore, $\kappa\in \C(X)$ (see Remark~\ref{obs_calibre_cofinito}).

Finally, Lemma~\ref{lema_calibres_cofinitos}(3) and the relations $\kappa>\omega$ and $\kappa\in \C(X)$ allow us to conclude that $$\CN \subseteq \{\lambda\in\CN : \lambda \leq \kappa\ \vee\ \cf(\lambda)>\omega\}\cup\{\lambda\in\CN : \cf(\lambda)<\kappa<\lambda \} \subseteq \C(X).$$
\end{proof}

As a consequence of Proposition~\ref{prop_cofinito_omega} we get a result for all countably infinite $T_1$ spaces. We need the following auxiliary lemma whose proof is evident.

\begin{lemma}\label{lema_calibres_topologias} Let $\tau$ and $\sigma$ be a pair of topologies on a set $X$. If $\sigma\subseteq \tau$, then $\C(X,\tau) \subseteq \C(X,\sigma)$, $\P(X,\tau) \subseteq \P(X,\sigma)$, $\WP(X,\tau) \subseteq \WP(X,\sigma)$, $\C^{*}(X,\tau) \subseteq \C^{*}(X,\sigma)$, $\P^{*}(X,\tau) \subseteq \P^{*}(X,\sigma)$ and $\WP^{*}(X,\tau) \subseteq \WP^{*}(X,\sigma)$.

\end{lemma}

\begin{theorem}\label{thm_calibres_espacios_numerables} If $X$ is a countably infinite $T_1$ space, then $\{\lambda \in\CN : \cf(\lambda)>\omega\} = \C(X) \subseteq \C^{*}(X) \subseteq \{\lambda \in\CN : \lambda>\omega\}$.

\end{theorem}

\begin{proof} On the one hand, since any $T_1$ topology extends the cofinite topology, Proposition~\ref{prop_cofinito_omega} and Lemma~\ref{lema_calibres_topologias} imply the inclusions $\C(X) \subseteq \{\lambda \in\CN : \cf(\lambda)>\omega\}$ and $\C^{*}(X) \subseteq \{\lambda \in\CN : \lambda>\omega\}$. On the other hand, since $X$ is a separable space, Lemma~\ref{lema_densidad} guarantees that $\{\lambda \in\CN : \cf(\lambda)>\omega\} \subseteq \CN$. Finally, the relation $\C(X) \subseteq \C^{*}(X)$ always holds.
\end{proof}

A natural question is whether in Theorem~\ref{thm_calibres_espacios_numerables} the equality $\C^{*}(X) = \{\lambda \in\CN : \lambda>\omega\}$ holds. We close this section with a theorem in which we show that this relation is independent of \textsf{ZFC}. Recall that the Continuum Hypothesis, \textsf{CH}, is the equality {\lq\lq}$\mathfrak{c}=\omega_1${\rq\rq}.

\begin{theorem}\label{thm_calibre*_independencia} The following statements hold.

\begin{enumerate}
\item If \textsf{CH} is true, then any countably infinite $T_1$ space $X$ satisfies the relation $\C^{*}(X) = \left\{\lambda \in\CN : \lambda >\omega\right\}$.

\item When $\aleph_\omega < \mathfrak{c}$, there exists a countably infinite metrizable space $X$ such that $\C^{*}(X) \neq \{\lambda \in\CN : \lambda>\omega\}$.
\end{enumerate}

\end{theorem}

\begin{proof} To verify item (1) let $X$ be a countably infinite $T_1$ space. An immediate observation is that $o(X) \leq \mathfrak{c}$ since $\tau_X$ is a family of subsets of $X$. Thus, under \textsf{CH} it is satisfied that $o(X) \leq \omega_1$. Finally, use the relations $\omega_1 \in \C(X) \subseteq \C^{*}(X)$ and $\{\lambda \in \CN : \lambda>o(X)\} \subseteq \C^{*}(X)$ to obtain that $\C^{*}(X) = \{\lambda \in\CN : \lambda>\omega \}$.

For item (2) let $X$ be the set $\omega$ equipped with the discrete topology. Use Lemma~\ref{lema_familia_casi_ajena} to find an almost disjoint family $\mathcal{A}$ such that $|\mathcal{A}|=\aleph_\omega$ and, for any $\mathcal{B} \in [\mathcal{A}]^{\aleph_\omega}$, there exist $A,B,C\in\mathcal{B}$ such that $A\cap B \cap C = \emptyset$. Under these circumstances, since $\mathcal{A}$ is a subset of $\tau_X^+$ of size $\aleph_\omega$, we deduce that $\aleph_\omega$ is not a precaliber* for $X$. Consequently, $\C^{*}(X) \neq \{\lambda \in\CN : \lambda>\omega\}$.
\end{proof}

\section{Chain conditions in topological sums}

For the purposes of this section, $\lambda$ will always represent a positive cardinal number, not necessarily infinite. The following result is a simple generalization of \cite[T.276, p.~311]{tkachuk2014}.

\begin{proposition}\label{preservacion_union} Let $X$ be a topological space and $\{X_\alpha : \alpha<\lambda\}$ be a family of subspaces of $X$. If $\cf(\kappa)>\lambda$ and $\kappa$ is a caliber (resp., precaliber; weak precaliber) for every $X_\alpha$, then $\kappa$ is a caliber (resp., precaliber; weak precaliber) for $\bigcup\{X_\alpha : \alpha<\lambda\}$.

\end{proposition}

\begin{proof} Set $Y:=\bigcup\{X_\alpha : \alpha<\lambda\}$ and let $\{U_\gamma : \gamma<\kappa\}$ be a subset of $\tau_Y^+$. For every $\gamma<\kappa$, define $\alpha(\gamma):=\min\left\{\alpha<\lambda : U_\gamma\cap X_\alpha\neq\emptyset\right\}$. Since $\cf(\kappa)>\lambda$, the function $\kappa \to \lambda$ given by $\gamma\mapsto \alpha(\gamma)$ has a fiber of size $\kappa$; i.e., there are $J\in [\kappa]^{\kappa}$ and $\alpha<\lambda$ such that $\alpha(\gamma)=\alpha$, whenever $\gamma\in J$. Finally, since $\{U_\gamma\cap X_\alpha : \gamma\in J\}$ is a subset of $\tau_{X_\alpha}^{+}$ and $X_\alpha$ has caliber $\kappa$, we get $I\in [J]^{\kappa}$ with $\bigcap\{U_\gamma\cap X_\alpha : \gamma\in I\}\neq\emptyset$; thus, $\{U_\gamma : \gamma\in I\}$ has non-empty intersection. The proof for precalibers and weak precalibers is analogous.
\end{proof}

\begin{proposition}\label{calibres_sumas_topologicas} If $\kappa$ is a cardinal number and $\{X_\alpha : \alpha<\lambda\}$ is a family of topological spaces, then $\bigoplus\{X_\alpha : \alpha<\lambda \}$ has caliber (resp., precaliber; weak precaliber) $\kappa$ if and only if $\cf(\kappa)>\lambda$ and each summand has caliber (resp., precaliber; weak precaliber) $\kappa$.

\end{proposition}

\begin{proof} The inverse implication is a consequence of Proposition~\ref{preservacion_union}. For the direct implication, suppose that $X:=\bigoplus\{X_\alpha : \alpha<\lambda \}$ has caliber (resp., precaliber; weak precaliber) $\kappa$. Since all $X_\alpha$ are open subspaces of $X$, then each summand has caliber (resp., precaliber; weak precaliber) $\kappa$ by Proposition~\ref{calibre_subconjunto_abierto}.

Finally, if $\cf(\kappa)\leq \lambda$, then $\cf(\kappa)$ is not a weak precaliber for $X$ since $\{X_\alpha : \alpha<\cf(\kappa)\}$ is a cellular family. Therefore, by Lemma~\ref{lema_cofinalidad}, $\kappa$ is not a weak precaliber for $X$ either, contrary to our hypotheses. Consequently, $\cf(\kappa)>\lambda$.
\end{proof}

We next look at the behaviour of chain* conditions under topological sums. We start with the following result which is consequence of Proposition~\ref{preservacion*_abiertos}.

\begin{proposition}\label{prop_calibre_baja_en_suma} If $\{X_\alpha : \alpha<\lambda\}$ is a family of topological spaces and $\kappa$ is a caliber* (resp., precaliber*; weak precaliber*) for $\bigoplus\{X_\alpha : \alpha<\lambda\}$, then $\kappa$ is a caliber* (resp., precaliber*; weak precaliber*) for $X_\alpha$, for every $\alpha<\lambda$.

\end{proposition}

The next result follows from Lemma~\ref{equiv_calibre_indices} and Proposition~\ref{preservacion_union}.

\begin{proposition}\label{preservacion*_union_regular} Let $X$ be a topological space and $\{X_\alpha : \alpha<\lambda\}$ be a familiy of subspaces of $X$. If $\kappa$ is a regular cardinal with $\kappa>\lambda$ and $\kappa$ is a caliber* (resp., precaliber*; weak precaliber*) for each $X_\alpha$, then $\kappa$ is a caliber* (resp., precaliber*; weak precaliber*) for $\bigcup\{X_\alpha : \alpha<\lambda\}$.

\end{proposition}

When $\kappa$ is a singular cardinal we have the following result.

\begin{theorem}\label{preservacion*_union} Let $\kappa$ be a singular cardinal with $\cf(\kappa)>\lambda$ and such that, for every regular cardinal $\cf(\kappa)\leq \mu<\kappa$, $\mu^\lambda=\mu$. Furthermore, let $X$ be a topological space and $\{X_\alpha : \alpha<\lambda\}$ be a family of subspaces of $X$. If $\kappa$ is a caliber* (resp., precaliber*; weak precaliber*) for every $X_\alpha$, then $\kappa$ is a caliber* (resp., precaliber*; weak precaliber*) for $\bigcup\{X_\alpha : \alpha<\lambda\}$.

\end{theorem}

\begin{proof} Set $Y:= \bigcup\{X_\alpha : \alpha<\lambda\}$ and fix a subset $\{U_\gamma : \gamma<\kappa\}$ of $\tau_Y^+$ enumerated without repetitions.

\medskip

\noindent {\bf Claim.} There is $\alpha<\lambda$ such that $\left|\left\{U_\gamma\cap X_\alpha : \gamma<\kappa\right\}\right|=\kappa$.

\medskip

Otherwise, for each $\alpha<\lambda$ the cardinal $\kappa_\alpha := \left|\left\{U_\gamma\cap X_\alpha : \gamma<\kappa\right\}\right|$ is strictly less than $\kappa$. So, $\cf(\kappa)>\lambda$ implies $\sup\{\kappa_\alpha : \alpha<\lambda\}<\kappa$ and thus, since $\kappa$ is a limit cardinal, there exists a regular $\mu<\kappa$ such that $\max\{\cf(\kappa), \sup\{\kappa_\alpha : \alpha<\lambda\}\} < \mu$.

Now, since the function $\{U_\gamma : \gamma<\kappa\} \to \{(U_\gamma \cap X_\alpha)_{\alpha<\lambda} : \gamma<\kappa\}$ given by $U_\gamma \mapsto (U_\gamma \cap X_\alpha)_{\alpha<\lambda}$ is one-to-one, and $\{(U_\gamma \cap X_\alpha)_{\alpha<\lambda} : \gamma<\kappa\}$ is a subset of the product $\prod_{\alpha<\lambda} \{U_\gamma \cap X_\alpha : \gamma<\kappa\}$, we deduce the relations \[\kappa \leq \left|\prod_{\alpha<\lambda} \{U_\gamma \cap X_\alpha : \gamma<\kappa\}\right| \leq  \sup\{\kappa_\alpha : \alpha<\lambda\}^{\lambda} \leq \mu^{\lambda}.\] Hence, since our hypotheses indicate that $\mu^{\lambda}=\mu$, we obtain the inequality $\kappa\leq \mu$, contradicting our assumption on $\mu$.

To finish our proof, if $\alpha<\lambda$ is as in our Claim, then Lemma~\ref{lema_equivalencia*} implies the existence of $J\in[\kappa]^{\kappa}$ with the following characteristics:
\begin{enumerate}
\item $\{U_\gamma\cap X_\alpha : \gamma\in J\}=\{U_\gamma\cap X_\alpha : \gamma<\kappa\}\setminus\{\emptyset\}$;
\item if $\gamma,\delta\in J$ are distinct, then $(U_\gamma\cap X_\alpha) \neq (U_\delta\cap X_\alpha)$; and
\item $\{U_\gamma\cap X_\alpha : \gamma\in J\}$ has non-empty intersection (resp., is centered; is linked).
\end{enumerate}
Consequently, $\left|\left\{U_\gamma : \gamma\in J\right\}\right|=\kappa$ and $\left\{U_\gamma : \gamma\in J\right\}$ has non-empty intersection (resp., is centered; is linked). Thus, $Y$ has caliber* (resp., precaliber*; weak precaliber*) $\kappa$.
\end{proof}

In the next result, the direct implication is a consequence of Proposition~\ref{prop_calibre_baja_en_suma}, and the converse  follows from Proposition~\ref{preservacion*_union_regular} and Theorem~\ref{preservacion*_union}.

\begin{theorem}\label{calibres*_sumas_topologicas} Let $\kappa$ be a cardinal number with $\cf(\kappa)>\lambda$ and such that, for every regular cardinal $\cf(\kappa)\leq \mu<\kappa$, $\mu^\lambda=\mu$. If $\{X_\alpha : \alpha<\lambda\}$ is a family of topological spaces, then $\bigoplus\{X_\alpha : \alpha<\lambda \}$ has caliber* (resp., precaliber*; weak precaliber*) $\kappa$ if and only if each summand has caliber* (resp., precaliber*; weak precaliber*) $\kappa$.

\end{theorem}

Recall that if $\mu$ is a regular cardinal with $\lambda<\mu$, then \textsf{GCH} implies that $\mu^{\lambda} = \mu$ (see \cite[Lemma~10.42, p.~34]{kunen1980}). Thus, we get (cf. Proposition~\ref{calibres_sumas_topologicas}):

\begin{theorem}\label{calibres*_sumas_topologicas_GCH} $[\mathsf{GCH}]$ If $\kappa$ is cardinal number with $\cf(\kappa)>\lambda$ and $\{X_\alpha : \alpha<\lambda\}$ is a family of topological spaces, then $\bigoplus\{X_\alpha : \alpha<\lambda \}$ has caliber* (resp., precaliber*; weak precaliber*) $\kappa$ if and only if each summand has caliber* (resp., precaliber*; weak precaliber*) $\kappa$.

\end{theorem}

With Theorem~\ref{calibres*_sumas_topologicas} available, we can obtain within \textsf{ZFC} a result similar to Theorems~\ref{thm_1} and \ref{thm_1.2} for uncountable cardinals with countable cofinality.

\begin{theorem}\label{thm_1_p2} If $\kappa$ is an uncountable cardinal with $\cf(\kappa)=\omega$, then there is a $T_1$ compact space $X$ such that $o(X)=\kappa = |X|$ and $\kappa \in \C^{*}(X)\setminus\C(X)$.

\end{theorem}

\begin{proof} Let $Y$ be the cofinite space of size $\kappa$, $Z$ the cofinite space of size $\omega$, and set $X:= Y\oplus Z$. Clearly, $X$ is compact, $T_1$ and $o(X)=\kappa = |X|$. Proposition~\ref{prop_cofinito_omega} shows that $\omega$ is not a caliber for $Z$, so that $\omega\not\in \C(X)$ by Proposition~\ref{calibre_subconjunto_abierto}. Consequently, $\kappa\not\in \C(X)$ by Lemma~\ref{lema_cofinalidad}. Finally, Propositions~\ref{prop_cofinito_omega} and \ref{prop_cofinito_>omega}, and Theorem~\ref{calibres*_sumas_topologicas}, imply that $\kappa\in \C^*(X)$.
\end{proof}

With an additional hypothesis, we can obtain a variant of Theorem~\ref{thm_1_p2} in which $\kappa\not\in \WP(X)$ is achieved.

\begin{theorem}\label{thm_1_p3} If $\kappa$ is an uncountable cardinal with $\kappa>\mathfrak{c}$ and $\cf(\kappa)=\omega$, then there is a $T_1$ locally compact space $X$ such that $o(X)=\kappa = |X|$ and $\kappa \in \C^{*}(X)\setminus\WP(X)$.

\end{theorem}

\begin{proof} Let $Y$ be the cofinite space of size $\kappa$, $Z$ the discrete space of size $\omega$, and set $X:= Y\oplus Z$. Evidently, $X$ is locally compact, $T_1$ and $o(X) = \kappa = |X|$. Now, since in Example~\ref{ejem_discreto} we proved that $\omega$ is not a weak precaliber for $Z$, Proposition~\ref{calibre_subconjunto_abierto} implies that $\omega\not\in\WP(X)$ and thus, $\kappa\not\in\WP(X)$ (see Lemma~\ref{lema_cofinalidad}). Finally, by Propositions~\ref{calibre*_topologia} and \ref{prop_cofinito_>omega}, and Theorem~\ref{calibres*_sumas_topologicas}, we conclude that $\kappa\in\C^*(X)$.
\end{proof}

Note that if $\omega\leq \kappa\leq\lambda$ is a cardinal number, then clearly $\kappa$ is not a weak precaliber* for $\bigoplus\{X_\alpha : \alpha<\lambda \}$, even if $\kappa$ is a caliber for every $X_\alpha$. In our following example we will show that the condition $\cf(\kappa)>\lambda$ is essential in Theorem~\ref{calibres*_sumas_topologicas_GCH}, even when $\kappa>\lambda$.

\begin{example} If $\cf(\lambda)>\omega$, then there is a family of compact Hausdorff spaces $\{X_\alpha : \alpha<\lambda\}$ such that, for every $\alpha<\lambda$, $o(X_\alpha)= 2^{\aleph_\lambda}=|X_\alpha|$, $\aleph_\lambda \in \C(X_\alpha)$, and $\aleph_\lambda \not \in \WP^{*}(\bigoplus\{X_\alpha : \alpha<\lambda\})$.

\end{example}

\begin{proof} Let $X_\alpha$ be the Cantor cube $D(2)^{\aleph_\lambda}$, for every $\alpha<\lambda$. Observe that Lemmas~\ref{lema_cubo_calibre} and \ref{lema_cubo_topologia} give that $\aleph_\lambda\in \C(X_\alpha)$ and $o(X_\alpha)=2^{\aleph_\lambda}=\abs{X_\alpha}$, while Lemma~\ref{lema_celular_no_WP} gives that $\aleph_\lambda \not \in \WP^{*}(\bigoplus\{X_\alpha : \alpha<\lambda\})$.
\end{proof}

We proved in Theorem~\ref{thm_1_p2} that if $X$ stands for the topological sum of the cofinite space of size $\aleph_\omega$ with the cofinite space of size $\omega$, then $\aleph_\omega \in \C^{*}(X)$. In this particular example, it is a consequence of Propositions~\ref{preservacion*_abiertos} and \ref{prop_cofinito_omega} that $\omega \not \in \C^{*}(X)$. For this reason, Theorem~\ref{thm_extension} implies the following corollary.

\begin{corollary}\label{cor2_thm_extension} There is a topological space $X$ and an extension $Y$ of $X$ such that $\aleph_\omega \in \C^{*}(X)$ and $\aleph_\omega \not\in \C^{*}(Y)$.

\end{corollary}

\section{Chain conditions in topological products}

There are several results about the preservation of the classical chain conditions for topological products  in the literature (see, for example, \cite{argtsa1982}, \cite{rios2022}, \cite{sanin1948} and \cite{shelah1977}). In this section we will present some preservation results related to chain* conditions.

\begin{lemma}\label{lema_funciones*} Let $f:X\to Y$ be a continuous surjective function. If $\kappa$ is a caliber* (resp., precaliber*; weak precaliber*) for $X$, then $\kappa$ is a caliber* (resp., precaliber*; weak precaliber*) for $Y$. 

\end{lemma}

\begin{proof} If $\{V_\alpha : \alpha<\kappa\}$ is a faithfully indexed subset of $\tau_Y^+$, then, since $\{f^{-1}[V_\alpha] : \alpha<\kappa\}$ is a subset of $\tau_X^+$ enumerated without repetitions, the hypothesis guarantees the existence of $J\in [\kappa]^{\kappa}$ such that $\{f^{-1}[V_\alpha] : \alpha\in J\}$ has non-empty intersection (resp., is centered; is linked). In this case, it is easy to see that $\{V_\alpha : \alpha\in J\}$ has non-empty intersection (resp., is centered; is linked).
\end{proof}

\begin{corollary}\label{cor_proyecciones} If $\{X_\alpha : \alpha<\lambda\}$ is a family of topological spaces and $\kappa$ is a caliber* (resp., precaliber*; weak precaliber*) for $\prod\{X_\alpha : \alpha<\lambda\}$, then, for every $\alpha<\lambda$, $\kappa$ is a caliber* (resp., precaliber*; weak precaliber*) for $X_\alpha$. 

\end{corollary}

In what follows we will show that chain conditions and chain* conditions {\lq\lq}truly{\rq\rq} coincide on Cantor cubes, i.e., if $\kappa$ is a caliber* for $D(2)^\lambda$ with $\kappa \leq 2^\lambda$, then $\kappa$ is a caliber for $D(2)^\lambda$.

\begin{theorem}\label{theo:cf}
Let $\kappa$ be a cardinal number with $\cf(\kappa)=\omega$ and $Z$ a topological space with $o(Z) \geq \kappa$. If $X \cong Y \times Z$, where $Y$ is an infinite Hausdorff space, then $X$ does not have weak precaliber* $\kappa$. 
\end{theorem}

\begin{proof} Since $o(Z) \geq \kappa$ and $Y$ is an infinite Hausdorff space, there are subsets $\{U_\alpha : \alpha <\kappa\} \subseteq \tau_Z^+$ and $\{V_n : n<\omega\} \subseteq \tau_Y^+$ such that $U_\alpha \neq U_\beta$ and $V_n \cap V_m = \emptyset$, whenever $\alpha<\beta<\kappa$ and $n<m<\omega$. Let $\{\kappa_n : n<\omega\}$ be a strictly increasing sequence of cardinal numbers such that $\kappa_0 := 0$ and $\sup\{\kappa_n : n<\omega\} = \kappa$.

For every $n < \omega$ define $\mathcal{W}_n := \{V_n \times U_\alpha : \kappa_n \leq \alpha < \kappa_{n+1}\}$. Then $\mathcal{W} := \bigcup_{n < \omega}\mathcal{W}_n$ is a family of $\kappa$ non-empty open subsets of $Y \times Z$. We claim that for every subset $\mathcal{W}' \subseteq \mathcal{W}$ of size $\kappa$, we can find two distinct elements of $\mathcal{W}'$ with empty intersection. Indeed, since $|\mathcal{W}'|=\kappa$, there are $n < m < \omega$ such that $\mathcal{W}'$ contains an element of $\mathcal{W}_n$ and an element of $\mathcal{W}_m$, say $V_n \times U_\alpha$ and $V_m \times U_\beta$. Finally, since $V_n \cap V_m = \emptyset$, then $(V_n \times U_\alpha) \cap (V_m \times U_\beta) = \emptyset$.
\end{proof}

For every infinite cardinal $\lambda$, $D(2)^\lambda \cong D(2)^\lambda \times D(2)^\lambda$, therefore we obtain: 

\begin{corollary}\label{cor_theo:cf} If $\kappa$ is a cardinal number with $\kappa \leq 2^\lambda$ and $\cf(\kappa)=\omega$, then $D(2)^\lambda$ does not have weak precaliber* $\kappa$. 
\end{corollary} 

\begin{corollary}\label{cor_condiciones_cadena_cubo} If $\lambda$ is a cardinal number, then:

\begin{enumerate}
\item $\C\left(D(2)^{\lambda}\right)=\P\left(D(2)^{\lambda}\right)=\WP\left(D(2)^{\lambda}\right)=\left\{\kappa\in\CN : \cf(\kappa)>\omega\right\}$.
\item \begin{align*} \C^*\left(D(2)^{\lambda}\right)&=\P^*\left(D(2)^{\lambda}\right)=\WP^*\left(D(2)^{\lambda}\right) \\
&=\left\{\kappa\in\CN : \kappa> 2^{\lambda}\right\}\cup \left\{\kappa\in\CN : \kappa\leq 2^{\lambda} \wedge \cf(\kappa)>\omega\right\}.
\end{align*}

\end{enumerate}

\end{corollary}

\begin{proof} To verify the first equalities, it is enough to observe that if $\kappa$ is a cardinal number with $\cf(\kappa)>\omega$, then Lemma~\ref{lema_cubo_calibre} implies that $\kappa$ is a caliber for $D(2)^{\lambda}$. On the other hand, since $D(2)^{\lambda}$ is an infinite Hausdorff space, no infinite cardinal with countable cofinality can be a weak precaliber for $D(2)^{\lambda}$ (see Lemmas~\ref{lema_cofinalidad} and~\ref{precalibre_debil_celularidad}).

For the second part, the relations $\left\{\kappa\in\CN : \cf(\kappa)>\omega\right\} = \C\left(D(2)^{\lambda}\right)\subseteq \C^*\left(D(2)^\lambda\right)$ together with Lemma~\ref{lema_cubo_topologia} and Proposition~\ref{calibre*_topologia} imply that $$\left\{\kappa\in\CN : \kappa > 2^{\lambda}\right\}\cup
\left\{\kappa\in\CN : \kappa \leq 2^{\lambda} \wedge \cf(\kappa)>\omega\right\} \subseteq \C^*(D(2)^\lambda).$$
Lastly, Corollary~\ref{cor_theo:cf} ensures that no $\kappa \leq 2^\lambda$ with countable cofinality is a weak precaliber* for $D(2)^\lambda$.
\end{proof}

We can extend Corollary~\ref{cor_condiciones_cadena_cubo} for certain classes of infinite Hausdorff spaces. In order to do this we first need to prove some auxiliary results regarding the behavior of the cardinal function $o$ in topological products. Let us consider the following question.

\begin{question} Let $X$ be a topological space and $\lambda$ be a cardinal number. Under what conditions on $X$ and $\lambda$ is it possible to express $o(X^\lambda)$ in terms of $\lambda$ and $\phi(X)$ for some topological cardinal function $\phi$? In particular, under what conditions on $X$ is $o(X^\lambda) = 2^{w(X)\cdot \lambda}$ (resp., $o(X^\lambda) = 2^{ |X|\cdot \lambda}$, $o(X^\lambda) = 2^{o(X)\cdot \lambda}$, $o(X^\lambda) = o(X)^\lambda$)?
\end{question}

\begin{lemma}\label{lema_finito} Let $\lambda$ be an infinite cardinal. If $Y:=  \prod_{\alpha < \lambda}Y_\alpha$, where $Y_\alpha$ is a Hausdorff (discrete) space with $1 < |Y_\alpha| < \omega$ for every $\alpha<\lambda$, then $o(Y) = 2^\lambda$.
\end{lemma}

\begin{proof} On the one hand, since $D(2)^\lambda$ embeds in $Y$ and $Y$ is $T_2$, $2^\lambda \leq o(Y)$ (see \cite[Theorem~3.1(b), p.~10]{hodel1984}). On the other hand, since each $Y_\alpha$ can be embedded as a subspace of the cube $D(2)^{\lambda}$, then $Y$ can be embedded as a subspace of the cube $\left(D(2)^\lambda\right)^{\lambda}\cong D(2)^\lambda$; consequently $o(Y) \leq o\left(D(2)^{\lambda}\right)=2^\lambda$ (see Lemma~\ref{lema_cubo_topologia}).
\end{proof}

From now on, the phrase {\lq\lq}$\phi$ is a topological cardinal function{\rq\rq} means that $\phi$ is a correspondence rule that assigns to each topological space $X$ an infinite cardinal number $\phi(X)$ such that, if $X$ is homeomorphic to $Y$, then $\phi(X) = \phi(Y)$.

Let $\mathcal{P}$ be a class of topological spaces and $\phi$ be a topological cardinal function. We will say that $(\mathcal{P},\phi)$ is a {\it finitely productive pair} if $\mathcal{P}$ is closed under finite products and, for each $X,Y\in \mathcal{P}$, $o(X) \leq 2^{\phi(X)}$ and $\phi(X\times Y) = \phi(X)\cdot\phi(Y)$. Similarly, $(\mathcal{P},\phi)$ will be called a {\it productive pair} if $\mathcal{P}$ is closed under arbitrary products and, for every $Y\in \mathcal{P}$ and $\{Y_\alpha : \alpha<\lambda\} \subseteq \mathcal{P}$, $o(Y) \leq 2^{\phi(Y)}$ and $\phi\left(\prod_{\alpha<\lambda} Y_\alpha\right) = \lambda\cdot\sup\left\{\phi(Y_\alpha) : \alpha<\lambda\right\}$.

These conventions will allow us to formulate with more simplicity the results that we will present below.

\begin{proposition}\label{prop_o(X)} Let $(\mathcal{P},\phi)$ be a finitely productive pair. If $X,Y \in \mathcal{P}$, $o(X)=2^{\phi(X)}$ and $o(Y)=2^{\phi(Y)}$, then $o(X\times Y) = o(X)\cdot o(Y)$, i.e., $o(X\times Y) = 2^{\phi(X)\cdot \phi(Y)}$.
 
\end{proposition}

\begin{proof} Clearly, $o(X \times Y) \leq 2^{\phi(X \times Y)} = 2^{\phi(X)\cdot\phi(Y)}$. On the other hand, the rule $\tau_X\times \tau_Y \to \tau_{X\times Y}$ defined by $(U,V)\mapsto U\times V$ determines a one-to-one function; hence, $2^{\phi(X) \cdot \phi(Y)} = 2^{\phi(X)} \cdot 2^{\phi(Y)} = o(X) \cdot o(Y) \leq o(X \times Y)$.
\end{proof}

\begin{corollary} Let $\phi \in \{w,nw\}$ (see \cite{hodel1984}). If $X$ and $Y$ are topological spaces such that $o(X)=2^{\phi(X)}$ and $o(Y)=2^{\phi(Y)}$, then $o(X\times Y) = o(X)\cdot o(Y) = 2^{\phi(X)\cdot \phi(Y)}$. Furthermore, the same conclusion is obtained when $X$ and $Y$ are infinite and $\phi = |\cdot|$.

\end{corollary}

Proposition~\ref{prop_o(X)} suggests the following natural question:

\begin{question}\label{Q_o(XxY)} Is it true that if $X$ satisfies $o(X) = 2^\kappa$ and $Y$ satisfies $o(Y) = 2^\lambda$, then $o(X \times Y)$ is of the form $2^\theta$ for some cardinal $\theta$?
\end{question}

Let $X$ and $Y$ be a pair of topological spaces. As already noted, the function $\tau_X\times \tau_Y \to \tau_{X\times Y}$ given by $(U,V)\mapsto U\times V$ is one-to-one, which implies that $o(X)\cdot o (Y) \leq o(X\times Y)$. On the other hand, if $W\in \tau_{X\times Y}$, for any $w\in W$ there exist $U_w\in \tau_X$ and $V_w \in \tau_Y$ such that $w\in U_w \times V_w \subseteq W$. Thus, the function $\tau_{X\times Y} \to P(\tau_X)\times P(\tau_Y)$ determined by $W \mapsto \left(\left\{U_w : w\in W\right\}, \left\{V_w : w\in W\right\}\right)$ is injective and therefore, if $o(X)\geq \omega$ or $o(Y)\geq \omega$, then $o(X\times Y) \leq 2^{o(X)\cdot o(Y)}$.

To sum up, if $X$ and $Y$ satisfy $o(X)\geq \omega$ or $o(Y)\geq \omega$, then $o(X)\cdot o(Y)\leq o (X\times Y) \leq 2^{o(X)\cdot o(Y)}$. In particular, if $o(X) = 2^\kappa$, $o(Y) = 2^\lambda$ and $\mu:=2^{\kappa\cdot\lambda}$, we have $\mu \leq o(X\times Y) \leq 2^{\mu}$. For this reason, under \textsf{GCH} the answer to Question~\ref{Q_o(XxY)} is affirmative since, with this hypothesis, $o(X\times Y)\in\{\mu, 2^{ \mu}\}$.

The following is a generalization of Proposition~\ref{prop_o(X)}.

\begin{proposition}\label{prop_prod} Let $(\mathcal{P},\phi)$ be a productive pair. If $\{X_\alpha : \alpha<\lambda\}\subseteq \mathcal{P}$ and for every $\alpha<\lambda$, $o(X_\alpha) = 2^{\phi(X_\alpha)}$, then $o(X) = \prod_{\alpha < \lambda}o(X_\alpha) = 2^{\lambda \cdot \kappa}$, where $X := \prod_{\alpha < \lambda}X_\alpha$ and $\kappa := \sup\{\phi(X_\alpha) : \alpha<\lambda\}$.
\end{proposition}

\begin{proof} As usual we denote the $\alpha^{\text{th}}$ projection by $\pi_\alpha : X \to X_\alpha$, for every $\alpha<\lambda$. We next define an injective function $f: \prod_{\alpha<\lambda} \left(\tau_{X_\alpha}\setminus\{X_\alpha\}\right) \to \tau_X$ by $$f(u) : = \bigcup_{\alpha < \lambda}\pi_\alpha^{-1}\left[u(\alpha)\right],$$ for every $u\in \prod_{\alpha<\lambda} \left(\tau_{X_\alpha}\setminus\{X_\alpha\}\right)$.

To verify that $f$ is indeed one-to-one, let $u,v\in \prod_{\alpha<\lambda} \left(\tau_{X_\alpha}\setminus\{X_\alpha\}\right)$ be distinct and fix $\beta<\alpha$ in such a way that $u(\beta)\neq v(\beta)$. Suppose, without loss of generality, that there exists $z\in u(\beta)\setminus v(\beta)$. For each $\alpha\in \lambda\setminus \{\beta\}$, the condition $v(\alpha) \neq X_\alpha$ implies that $X_\alpha\setminus v(\alpha)\neq \emptyset$. Take any $x_\alpha \in X_\alpha\setminus v(\alpha)$ and consider the function $$x:= \left\{(\beta,z)\right\} \cup \left\{(\alpha,x_\alpha) : \alpha\in \lambda\setminus \{\beta\}\right\}\in X.$$ It is not difficult to see that $x$ belongs to $f(u)$ but is not an element of $f(v)$. Thus $f$ is injective and therefore, $\prod_{\alpha < \lambda}o(X_\alpha) \leq o(X)$.

Consequently, a routine cardinal arithmetic argument shows that \begin{align*} 2^{\lambda \cdot \kappa} &= 2^{\lambda \cdot \sup\{\phi(X_\alpha) : \alpha<\lambda\}} = 2^{\Sigma_{\alpha < \lambda}\phi(X_\alpha)} = \prod_{\alpha < \lambda}2^{\phi(X_\alpha)} \\
&= \prod_{\alpha < \lambda}o(X_\alpha) \leq o(X) \leq 2^{\phi(X)} = 2^{\lambda \cdot \sup\{\phi(X_\alpha) : \alpha<\lambda\}} = 2^{\lambda \cdot \kappa}.
\end{align*} In conclusion, $o(X) = \prod_{\alpha < \lambda}o(X_\alpha) = 2^{\lambda \cdot \kappa}$.
\end{proof}

\begin{corollary} Let $(\mathcal{P},\phi)$ be a productive pair. If $X\in \mathcal{P}$ and $o(X) = 2^{\phi(X)}$, then $o(X^\lambda) = o(X)^{\lambda} = 2^{\lambda \cdot \phi(X)} = 2^{\phi(X^\lambda)}$.
\end{corollary}

It is known that if there are no inaccessible cardinals and \textsf{GCH} is true, then any infinite Hausdorff space satisfies $o(X) = 2^\kappa$ for some cardinal $\kappa$ (see \cite[Theorem~ 13.1, p.~48]{hodel1984}). This observation combined with what we have discussed so far suggests the following question:

\begin{question}\label{Q_2^phi} If $X$ is a topological space and $\phi$ is a cardinal function, under what conditions is the equality $o(X)=2^{\phi(X)}$ satisfied?

\end{question}

For example, regarding the cardinal function $\phi=w$, in \cite[Theorem~8.1(e), p.~32]{hodel1984} it is shown that, if $X$ is infinite and metrizable, then $o(X) = 2^{w(X)}$. On the other hand, in the realm of compact Hausdorff spaces it is not always satisfied that $o(X) = 2^{w(X)}$. To illustrate this fact notice that, if $X$ is the Alexandroff-Urysohn double arrow (see \cite[Example~14.4, p.~51]{hodel1984}), then $o(X) = \mathfrak{c} = w(X)$. However, in certain subclasses of compact $T_2$ spaces we can obtain the relation $o(X) = 2^{w(X)}$.

A space is {\it extremely disconnected} if every open subset of it has an open closure. In \cite[Corollary~1, p.~608]{balfra1982} it is proved that extremely disconnected compact spaces satisfy the equality $|X| = 2^{w(X)}$. Consequently, since $|X|\leq o(X) \leq 2^{w(X)}$ for $T_1$-spaces, we conclude that $o(X)=2^{w(X)}$. For example, any Stone space of a complete Boolean algebra verifies this relationship (see \cite[Theorem~(d), p.~448]{porwoo1988}).

On the other hand, we say that a compact space is {\it polyadic} if it is a continuous image of a power of the one-point compactification of an infinite discrete space. If $X$ is a polyadic compact space, then $X$ admits a discrete subspace of cardinality $w(X)$ (see \cite[Corollary~1, p.~15]{gerlits1978}). Thus, since for each discrete subspace $D\subseteq X$ it is satisfied that $2^{|D|}\leq o(X)$, we deduce the equality $o(X)=2^{w(X)} $.

The class of polyadic compact spaces is broad. Indeed, observe that if $\kappa\geq\omega$, then the function $f:\alpha D(\kappa) \to D(2)$ determined by $f(0)=0$ and $f \left[\alpha D(\kappa)\setminus\{0\}\right]\subseteq \{1\}$ is continuous and surjective. For this reason, all dyadic compacts (i.e., the continuous images of Cantor cubes) are polyadic compact spaces; in particular, all compact topological groups belong to this class (see \cite{shakhmatov1994}).

Back to our original goal, with this background we are better positioned to generalize Corollary~\ref{cor_condiciones_cadena_cubo}.

\begin{lemma}\label{lema_producto} Let $\kappa$ be a cardinal number with $\cf(\kappa)=\omega$ and $\{X_\alpha : \alpha < \lambda\}$ be a family of topological spaces each containing more than one point. If there exists $J \subseteq \lambda$ such that $\prod_{\alpha \in J}X_\alpha$ is Hausdorff and infinite, and $o\left(\prod_{\alpha \in \lambda \setminus J}X_ \alpha\right) \geq \kappa$, then $\prod_{\alpha < \lambda}X_\alpha$ does not have weak precaliber* $\kappa$.
\end{lemma}

\begin{proof} It is enough to define $Y := \prod_{\alpha \in J}X_\alpha$ and $Z := \prod_{\alpha \in \lambda \setminus J}X_\alpha$, and apply Theorem~\ref{theo:cf}.
\end{proof}

\begin{lemma}\label{lema_Hausdorff} Let $(\mathcal{P},\phi)$ be a productive pair, $\kappa$ be a cardinal number with $\cf(\kappa)=\omega$ and $\{X_\alpha : \alpha < \lambda\}\subseteq \mathcal{P}$ be a family of Hausdorff spaces such that $o(X)\geq \kappa$, where $X:=\prod_{\alpha < \lambda}X_\alpha$. Under these hypotheses, if $\lambda\geq\omega$ and for each $\alpha<\lambda$ it is satisfied that $|X_\alpha|>1$ and $o(X_\alpha) = 2^{\phi( X_\alpha)}$, then $X$ does not have weak precaliber* $\kappa$. 

\end{lemma}

\begin{proof} Firstly, if every $X_\alpha$ is finite, then we take a countable $J\subseteq \lambda$ with $|\lambda \setminus J| \geq \omega$. It turns out that $o\left(\prod_{\alpha \in \lambda \setminus J}X_\alpha\right) = o(X)$ (see Lemma~\ref{lema_finito}) and thus, Lemma~\ref{lema_producto} implies that $X$ does not have weak precaliber* $\kappa$.

Otherwise, there exists $\beta<\lambda$ such that $X_{\beta}$ is infinite. If $o\left(\prod_{\alpha \in \lambda \setminus \{\beta\}}X_\alpha\right) = o(X)$, then Lemma~\ref{lema_producto} guarantees that $X$ does not have weak precaliber* $\kappa$. On the other hand, if $o\left(\prod_{\alpha \in \lambda \setminus \{\beta\}}X_\alpha\right) < o(X)$, then Proposition~\ref{prop_prod} gives $o(X) = o\left(\prod_{\alpha \in \lambda \setminus \{\beta\}}X_\alpha\right) \cdot o\left(X_{\beta}\right)$, and hence the equality $o\left(X_{\beta}\right)=o(X)$. Finally, we use Lemma~\ref{lema_producto} to conclude that $X$ does not have weak precaliber* $\kappa$.
\end{proof}

\begin{corollary}\label{cor_potencia} Let $(\mathcal{P},\phi)$ be a productive pair, $\kappa$ be a cardinal number with $\cf(\kappa)=\omega$ and $X\in \mathcal{P}$ be a Hausdorff space such that $o(X^\lambda)\geq \kappa$. Under these hypotheses, if $\lambda\geq\omega$, $|X|>1$ and $o(X) = 2^{\phi(X)}$, then $X^\lambda$ does not weak precaliber* $\kappa$.

\end{corollary}

In \cite{shelah1977} it was shown that, if $\kappa$ is an infinite cardinal, $X$ is a topological space with caliber $\kappa$, $\cf(\kappa)>\omega$ and $\lambda$ is a cardinal number, then the power $X^{\lambda}$ also has caliber $\kappa$. With this result available, we can prove the following generalization of Corollary~\ref{cor_condiciones_cadena_cubo}:

\begin{corollary} Let $(\mathcal{P},\phi)$ be a productive pair, $\lambda$ be an infinite cardinal and $X\in \mathcal{P}$ be a Hausdorff space with $|X|>1$ and $o(X) = 2^{\phi(X)}$. If $\mu := o\left(X^\lambda\right)$, then for every $\mathsf{A}\in \{\C,\P,\WP\}$ the following statements hold:

\begin{enumerate}
\item $\mathsf{A}\left(X^{\lambda}\right)=\left\{\kappa\in\mathsf{A}\left(X\right) : \cf(\kappa)>\omega\right\}$.
\item $\mathsf{A}^*\left(X^{\lambda}\right) = \left\{\kappa\in\CN : \kappa> \mu\right\}\cup \left\{\kappa\in\mathsf{A}(X) : \kappa\leq \mu \wedge \cf(\kappa)>\omega\right\}$.

\end{enumerate}

\end{corollary}

\begin{proof} Since the proof is similar for $\C$, $\P$ and $\WP$, we will only expose the details for $\C$. First, the result of \cite{shelah1977} implies that $\left\{\kappa\in\C(X) : \cf(\kappa)>\omega\right\} \subseteq \C\left(X^{\lambda} \right)$. On the other hand, if $\kappa\in \C\left(X^{\lambda} \right)$, then clearly $\kappa\in \C(X)$, while $\cf(\kappa)>\omega$ follows from the fact that $X^{\lambda}$ is an infinite Hausdorff space (see Lemmas~\ref{lema_cofinalidad} and~\ref{precalibre_debil_celularidad}).

For item (2) notice that the relations $\left\{\kappa\in\C(X) : \cf(\kappa)>\omega\right\} = \C\left(X^{\lambda}\right)\subseteq \C^*\left(X^\lambda\right)$ and Proposition~\ref{calibre*_topologia} guarantee that $$\left\{\kappa\in\CN : \kappa > \mu\right\}\cup
\left\{\kappa\in\C(X) : \kappa \leq \mu \wedge \cf(\kappa)>\omega\right\} \subseteq \C^*(X^\lambda).$$
Finally, Corollary~\ref{cor_potencia} confirms that there is no $\kappa \leq \mu$ with $\cf(\kappa)=\omega$ such that $\kappa$ is a caliber* for $X^\lambda$.
\end{proof}

Since the cardinal functions $w$ and $nw$ satisfy the conditions of the previous result (see \cite{juhasz1980}), we obtain our main result regarding topological powers:

\begin{corollary}\label{cor_potencia_H} If $\lambda$ is an infinite cardinal, $X$ is a Hausdorff space with $|X|>1$, $\phi\in \{w,nw\}$, $o(X) = 2^{\phi(X)}$ and $\mu := o\left(X^\lambda\right)$, then for every $\mathsf{A}\in \{\C,\P,\WP\}$ the following statements hold:

\begin{enumerate}
\item $\mathsf{A}\left(X^{\lambda}\right)=\left\{\kappa\in\mathsf{A}\left(X\right) : \cf(\kappa)>\omega\right\}$.
\item $\mathsf{A}^*\left(X^{\lambda}\right) = \left\{\kappa\in\CN : \kappa> \mu\right\}\cup \left\{\kappa\in\mathsf{A}(X) : \kappa\leq \mu \wedge \cf(\kappa)>\omega\right\}$.

\end{enumerate}

\end{corollary}

The authors thank the reviewer for many careful and detailed comments that greatly improved the quality of this article.

\end{document}